\documentclass[11pt]{amsart}
\usepackage{graphicx}
\usepackage{amssymb,amsmath}
\usepackage{xy}
\xyoption{all}
\xyoption{poly}

\newtheorem{theorem}{Theorem}[section]
\newtheorem{corollary}[theorem]{Corollary}
\newtheorem{proposition}[theorem]{Proposition}

\newtheorem{lemma}[theorem]{Lemma}
\newtheorem{definition}[theorem]{Definition}

\newcommand{\cal}{\mathcal}
\def\qed{\hfill \mbox{$\Box$} \vspace{3mm}}

\def\R{{\mathcal{R}}}

\def\FS{{\mathfrak{S}}}

\def\BC{{\mathbf{C}}}

\def\BN{{\mathbb{N}}}
\def\BC{{\mathbb{C}}}
\def\BR{{\mathbb{R}}}
\def\BZ{{\mathbb{Z}}}

\def\Abs{\operatorname{Abs}\nolimits}

\def\Id{\operatorname{Id}\nolimits}

\def\GL{\operatorname{GL}\nolimits}

\def\Red{\operatorname{Red}\nolimits}

\def\res{\operatorname{res}\nolimits}

\def\Diag{\operatorname{Diag}\nolimits}

\def\reg{\operatorname{reg}\nolimits}

\def\codim{\operatorname{codim}\nolimits}
\def\disc{\operatorname{disc}\nolimits}
\def\triv{\operatorname{triv}\nolimits}

\def\ie{{\em i.e.}}

\title{Non-crossing partitions of type $(e,e,r)$}
\author{David Bessis \& Ruth Corran}
\address{DMA, \'Ecole normale sup\'erieure, 
45 rue d'Ulm, 75230 Paris cedex 05, France}
\email{david dot bessis at ens dot fr}
\address{Institut Henri Poincar\'e, 11 rue Pierre et Marie Curie,
75231 Paris cedex 05, France}
\email{corran at ihp dot jussieu dot fr}

\begin{document}

\begin{abstract}
We investigate a new lattice of generalised non-crossing partitions,
constructed using the geometry of the complex reflection group
$G(e,e,r)$. For the particular case $e=2$ (resp. $r=2$), our lattice coincides
with the lattice of simple elements for the type $D_n$ (resp. $I_2(e)$)
dual braid monoid. 
Using this lattice, we construct a Garside structure for 
the braid group $B(e,e,r)$. As a corollary, one may solve 
the word and conjugacy problems in this group.
\end{abstract}

\maketitle


\section*{Introduction}

The object of this article is the study
of the braid group associated with the complex reflection group $G(e,e,r)$,
via a new presentation by generators and relations. This new presentation
was introduced in our earlier unpublished
work \cite{bc}, where we proved that
it defines a {\em Garside monoid}, which has many consequences:
the monoid embeds in the braid group, the braid
group admits nice normal forms, from which one may describe
solutions to the word and conjugacy problems. However, our initial proofs
in \cite{bc} involved computer assisted case-by-case analysis.
While completing
\cite{bc}, we realised that most computations could be avoided by using
a suitable notion of {\em non-crossing partitions}.
The non-crossing partitions of type $(e,e,r)$ form a lattice $NCP(e,e,r)$,
with remarkable numerical properties, and sharing many features with
the lattice of classical non-crossing partitions (which we denote
by $NCP(1,1,r)$, and which corresponds to the symmetric group).
The rich combinatorics of this
lattice reflect geometric properties of $G(e,e,r)$.
We use these combinatorics to obtain
precise structural information about the associated braid group $B(e,e,r)$,
but they might also be meaningful in other
areas of mathematics where non-crossing partitions have recently appeared,
such as cluster algebras and free probabilities.

\vspace{3mm}

The classification of finite complex
reflection groups was obtained fifty years ago by Shephard and Todd,
\cite{shto}.
The problem essentially amounts to classifying the 
irreducible groups, whose list consists of:
\begin{itemize}
\item an infinite family
$G(de,e,r)$, where $d,e,r$ are
arbitrary positive integral parameters;
\item $34$ exceptions, labelled $G_{4},\dots,G_{37}$.
\end{itemize}
The infinite family includes the four infinite families of finite Coxeter
groups: $G(1,1,r)\simeq W(A_{r-1})$, $G(2,1,r)\simeq W(B_r)$,
$G(2,2,r)\simeq W(D_r)$ and $G(e,e,2)\simeq W(I_2(e))$. For all other
values of the parameters, $G(de,e,r)$ is an irreducible monomial complex
reflection group of rank $r$, with no real structure.

A general problem is to generalise to complex reflection groups as much
as possible from the theory of Weyl groups and Coxeter groups.
One of the motivations is the \emph{Spetses} program of
Brou\'e-Malle-Michel (see \cite{bmm}),
which formally extends certain representation-theoretic
aspects of reductive algebraic groups, as if certain complex reflection
groups were the ``Weyl groups'' of more general structures
(the \emph{Spetses}).
Not all complex reflection groups give rise to {Spetses}. In
the infinite family, one considers, in addition to the real groups,
the complex subfamilies $G(d,1,r)$ and $G(e,e,r)$. Among recent developments
is the description by Brou\'e-Kim, \cite{brkim}, of the ``Lusztig 
families of characters'' of $G(d,1,r)$ and $G(e,e,r)$.
These representation-theoretic aspects involve generalised braid groups
and Hecke algebras, which were introduced by Brou\'e-Malle-Rouquier.
Following \cite{bmr}, one defines the \emph{braid group} $B(G)$
attached to a complex reflection group $G$ 
as the fundamental group of the space of regular orbits.

When $G$ is real, the braid group $B(G)$ is well understood thanks to
Brieskorn's presentation theorem and the subsequent
structural study by
Deligne and Brieskorn-Saito (\cite{brieskorn}, \cite{deligne},
\cite{brieskornsaito}).
The braid group of $G(d,1,r)$ is isomorphic to the type $B_r$
Artin group, hence is also well understood. 
Our object of interest is the braid group of $G(e,e,r)$.
Note that this subseries contains the $D$-type and $I_2$-type Coxeter
series.

Using an induction based on fibration arguments, a presentation for the
braid group of $G(e,e,r)$ was computed in \cite{bmr}. This presentation shares
with those of
finite type Artin groups (= braid groups of finite Coxeter groups)
the following features:
there are $r$ generators, which correspond to \emph{braid reflections}
(also called \emph{generators of the monodromy}), relations are positive and
homogeneous, and addition of quadratic relations yields a presentation
of the reflection group. But when both $e\geq 3$ and $n\geq 3$ the
monoid defined by this presentation fails to embed in the braid group
(see \cite{corranphd}, p. 122).
By contrast with finite type Artin presentations, Brou\'e-Malle-Rouquier's
presentation is not {\em Garside}, in the sense of Dehornoy-Paris,
\cite{dehpa}.
The techniques of Deligne and Brieskorn-Saito are not applicable
to this presentation.

In \cite{bkl}, \cite{bdm} and \cite{dualmonoid}, new Garside presentations
for finite type Artin group were introduced. Geometrically, the generating
sets naturally arise when ``looking'' at the reflection arrangement
from an eigenvector for the Coxeter element. In $G(e,e,r)$, maximal
order regular elements (in the sense of \cite{springer}) can be seen
as analogues of Coxeter elements. Our new presentation for $B(e,e,r)$
is analog of the presentation from \cite{dualmonoid}. The generators
are {\em local generators} associated with eigenvectors
of such generalised Coxeter elements (as in \emph{loc. cit.}, section 4).
In particular, when restricted to Coxeter types $D_r$ and $I_2(e)$,
our monoid coincides with the dual braid monoid.

The lattice $NCP(e,e,r)$ coincides with the lattice of \emph{simple elements}
of our monoid. Its numerology (cardinality, M\"obius number, Zeta function)
may be described in terms of the reflection degrees of $G(e,e,r)$,
following the pattern observed by Chapoton for the Coxeter types
(\cite{chapoton}). In particular, the cardinality of $NCP(e,e,r)$ is a 
\emph{generalised Catalan number}.
When $G$ is a Weyl group, the numerical invariants
of the dual braid monoid are connected with those of cluster algebras
(see again \cite{chapoton}).
This is another evidence suggesting to look for analogs of
representation-theoretic objects, where the Weyl group is replaced
by a spetsial complex reflection group.
Another area of mathematics where non-crossing partitions (associated
with $G(1,1,r)$ and $G(2,1,r)$) have recently appeared is
the theory of free probabilities (\cite{bgn}) -- is it possible
to use $NCP(e,e,r)$ to construct a free probability theory?

There are two possible ways of proving that a certain monoid is a
Garside monoid. One possibility is to use the notion
of {\em complete presentations} (see \cite{dehornoy});
this purely word-theoretical
approach, used in our earlier work \cite{bc}, is applicable in quite general 
situations, but does not provide intuitive explanation of why the result
is true. The second approach, used for example in \cite{bdm}, is to
find a good interpretation of the poset of simple elements, which allows
a direct proof of the lattice property. This second approach is used
here.
Our strategy resembles the one used in \cite{bdm} 
and \cite{dualmonoid}.

\smallskip

The first section contains the construction of
our lattice $NCP(e,e,r)$, and a list of its basic properties; 
it may be read independently from the rest of the article.
Given any finite subset $x\subset \BC$, one may define a 
natural notion of \emph{non-crossing partition of $x$}.
These partitions are naturally ordered by refinement.
The non-crossing partitions of type $(e,e,r)$ are particular non-crossing
partitions of the set $\mu_{e(r-1)}\cup \{0\}$ (a regular $e(r-1)$-gon,
together with its center $0$), which are symmetric or asymmetric
in a particular way.
Even though the lattice $NCP(e,e,r)$
is related to the geometry of $G(e,e,r)$,
its description is purely combinatorial.
The case $e=2$ gives a notion of non-crossing partitions of type $D_r$,
which does not
coincide with initial definition by Reiner (\cite{reiner}), but
is equivalent to the independent new 
definition by Athanasiadis-Reiner (\cite{ar})
and provides a geometric interpretation for the simple elements of
the dual monoid of type $D_r$. When $r=2$, we also obtain a
geometric interpretation for the dual monoid of type $I_2(e)$.
As expected,
there is an analog of Kreweras complement operation:
given $u\preccurlyeq v$ in $NCP(e,e,r)$, we define a new element
$u\backslash v$.

In the second section, we construct natural maps from $NCP(e,e,r)$ to
$B(e,e,r)$ and $G(e,e,r)$. 
The construction relies on the choice of a generalised Coxeter element
$c\in G(e,e,r)$.
The map $u\mapsto b_u$ satisfies the property
$b_u b_{u\backslash v} = b_v$ for all $u\preccurlyeq v$ in $NCP(e,e,r)$
(Proposition \ref{mappingb}). This provides an \emph{a posteriori}
justification for using a ``left-quotient'' notation for the complement
operation.

In section \ref{section3}, we have a closer look at minimal non-trivial
partitions and their images in $B(e,e,r)$ and $G(e,e,r)$. 
This prepares for section \ref{section4}, which contains a crucial
step in our construction.
Consider the set $T$ of all reflections in $G(e,e,r)$.
There is a notion of \emph{reduced $T$-decompositions} of an element
$g\in G$. Define a relation $\preccurlyeq_T$ on $G(e,e,r)$, by setting
$g\preccurlyeq_T h$ whenever $g$ is the product of an initial segment of 
a reduced $T$-decomposition of $h$.
We prove that the map from $NCP(e,e,r)$ to $G(e,e,r)$ lands in
$P_G:=\{g\in G(e,e,r) | g\preccurlyeq c\}$.
Moreover, it induces a poset isomorphism $(NCP(e,e,r),\preccurlyeq)
\simeq (P_G,\preccurlyeq_T)$.
In particular, $(P_G,\preccurlyeq_T)$ is a lattice,
generalising what happens in Coxeter groups.

Once this lattice property is proved, the Garside structure for $B(e,e,r)$
is obtained very easily (section
\ref{section5}). Besides the many general properties of
Garside structures (solution to the word and conjugacy problems, finite
classifying spaces, etc...) which we do not detail, we mention
another application: the computation of 
certain centralisers, following a conjecture
stated in \cite{bdm}.

Another application is the description of an explicit presentation.
In section \ref{section7}, we give such a presentation; generators
correspond to minimal non-discrete partitions, and relations are of length
$2$. This follows the pattern observed in \cite{dualmonoid} for
real reflection groups. Here again, the fact that it is enough
to consider relations of length $2$ is a consequence of the
transitivity of the classical braid group Hurwitz action on the
set of reduced $T$-decompositions of $c$ (section \ref{section6}).

Section \ref{section8} deals with enumerative aspects of $NCP(e,e,r)$.

In the final section, we return to the study of Brou\'e-Malle-Rouquier's
presentation. Not only does t not satisfy the embedding
property, but it is impossible to add a finite number of relations
(without changing the generating set)
to obtain a presentation with the embedding property.

\smallskip

For real reflection groups, both the classical ``Coxeter-Artin'' viewpoint and
the dual approach provide Garside structures. Here only the dual approach
is applicable. The main features of our construction actually generalise
to all well-generated complex reflection groups (see \cite{hurwitz},
and its sequel in preparation).

\section*{Notations and terminology}

Throughout this article, $e$ and $n$ are fixed positive integers.
Rather than considering $G(e,e,n)$, it simplifies notations to 
work with $G(e,e,n+1)$.

We denote by $|X|$ the cardinal of a set $X$.

If $d$ is a positive integer, we denote by $\mu_d$ the group
of complex $d$-th roots of unity. We denote by $\zeta_d$ the
generator $\exp(2i\pi/d)$ of $\mu_d$.

When $X$ is a finite set, a cyclic ordering on $X$ consists
of a subset of $X^3$ (whose elements are called \emph{direct triples})
subject to certain (obvious) axioms.
E.g., the set $\mu_d$ comes equipped with a standard counterclocwise
cyclic ordering.
Any cyclic ordering on a finite set $X$ is equivalent, via a bijection
$X\simeq \mu_{|X|}$, to this standard example (we may use this property
as a definition, and avoid stating explicitly the axioms).

When $X$ is endowed with a cyclic ordering and $x_1,x_2\in X$,
we denote by $\left<x_1,x_2\right>$ the set
of $x\in X$ such that $(x_1,x,x_2)$ is
direct. We use the similar notations $\left> x_1,x_2\right>$,
$\left<x_1,x_2\right<$ and
$\left>x_1,x_2\right<$ to exclude one or both of the endpoints.

In \cite{bdm} and \cite{dualmonoid}, the symbol $\prec$ is used
to denote various order relations. The symbol $\preccurlyeq$ is used
here for the same purposes. When we write $u\prec v$, we mean
$u\preccurlyeq v$ and $u\neq v$.

\section{Generalised non-crossing partitions}
\label{section1}

Our goal here is to construct a suitable theory of non-crossing partitions
of type $(e,e,n+1)$.

While our construction is motivated by geometric considerations
on complex reflection groups, it is possible to give purely combinatorial
definitions. Formally, the reader solely interested in
combinatorics of $NCP(e,e,n+1)$ can read this section without any knowledge
of complex reflection groups.

\subsection{Set-theoretical partitions}
Let $x$ be a set.
A {\em (set-theoretical) partition} of $x$ is an unordered family
$u=(a)_{a\in u}$ of disjoint non-empty subsets
of $x$ (the ``parts'') such that $x=\bigcup_{a\in u}a$.
Partitions are in natural
bijection with abstract equivalence relations on $x$.

The set $P_x$ of partitions of $x$ is endowed with the
following order relation: 
\begin{definition}
Let $u,v$ be partitions of $x$. We say
that ``$u$ refines $v$'', and write $u\preccurlyeq v$,
if and only if $$\forall a\in u,\exists
b\in v, a \subseteq b.$$
\end{definition}

We recall that a \emph{lattice} is a poset $(P,\leq)$ such that, for any pair
of elements $u,v\in P$, there exists elements $u\wedge v$
(read ``$u$ meet $v$'') and $u\vee v$ (read ``$u$ join $v$'')
in $P$,
such that $u\wedge v \leq u, u \wedge v \leq v$, and
$\forall w\in P, (w \leq u) \; \text{and} \; (w \leq v ) \Rightarrow
w \leq u\wedge v$, and symmetrically for $u\vee v$.

Clearly, $(P_x,\preccurlyeq)$ is a lattice.
For any $u,v\in P_x$, $u\wedge v$ is the partition whose parts
are the non-empty $a\cap b$, for $a\in u$ and $b\in v$.
The element $u\vee v$ is better understood in terms of equivalence
relations: it is the smallest equivalence relation whose graph contains
the graphs of $u$ and $v$.
The minimal element of $P_x$ is the {\em discrete} partition ``$\disc$''
whose parts
are singletons. The maximal element of $P_x$ is the {\em trivial}
partition ``$\triv$'' with the unique part $x$.

In the sequel, we will consider several subsets of $P_x$. They will
always be endowed with the order relation obtained by restricting
$\preccurlyeq$. Regardless of the subset of $P_x$, we will always denote this
restriction by $\preccurlyeq$.
Our point will be to prove that the considered subsets are lattices. In all
cases, it will be enough to rely on the following easy
lemma to
prove the lattice property.

\begin{lemma}
\label{triviallemma}
Let $(P,\leq)$ be a finite lattice. Let $Q\subseteq P$, endowed with the
restricted order.
Denote by $m$ the maximal element of $P$. Denote by $\wedge_P$ the
meet operation in $P$.
Assume that
$$m \in Q \qquad \text{and} \qquad
\forall u,v\in Q, u\wedge_P v \in Q.$$
Then $(Q,\leq)$ is a lattice, with meet operation $\wedge_Q = \wedge_P$, and
join operation $\vee_Q$ satisfying
$$\forall u, v\in Q, u\vee_Q v = \bigwedge_{
\begin{array}{c}
w\in Q \\
u\leq w,v\leq w
\end{array}} w.$$
\end{lemma}

In such a context, it is convenient to use the notation $\wedge$ to
refer to both $\wedge_P$ and $\wedge_Q$. 

\subsection{Non-crossing partitions of a configuration of points}

Let $x$ be a finite subset of $\BC$. 
A partition $\lambda$ of $x$ is said to be non-crossing if it satisfies
the following condition:
\begin{itemize}
\item[] \emph{For all $\nu,\nu'\in\lambda$, if the convex hulls of $\nu$ and
$\nu'$ have a non-empty intersection, then $\nu=\nu'$.}
\end{itemize}
We denote by $NCP_x$ the set of non-crossing partitions of $x$.

\begin{lemma}
\label{generalncp}
The poset $(NCP_x,\preccurlyeq)$ is a lattice.
\end{lemma}

\begin{proof}
By Lemma \ref{triviallemma}, it is enough to check that the trivial
partition is in $NCP_x$ and that
the meet of two non-crossing partitions is non-crossing. Both follow
immediately from the definition.
\end{proof}

\subsection{Non-crossing partitions of type $(1,1,n)$}

The notation $NCP(1,1,n)$ is just another name for the lattice of non-crossing
partitions of a regular $n$-gon; this very classical object is related to
the geometry of the symmetric group $\FS_n$, aka $W(A_{n-1})$ or $G(1,1,n)$:

\begin{definition}
For any positive integer $n$, we set $NCP(1,1,n):=NCP_{\mu_n}$.
\end{definition}

\begin{figure}[ht]
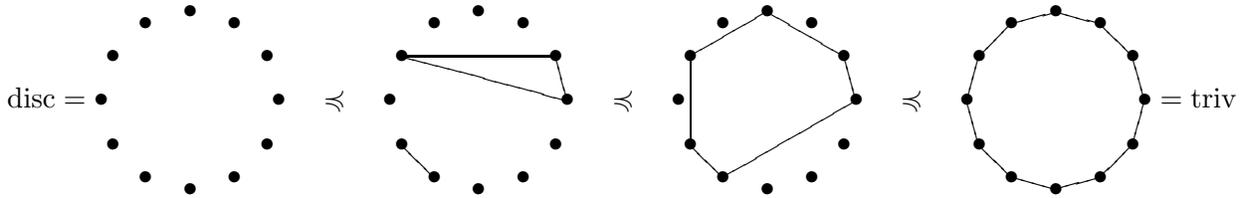

$$
\disc = \xy/r2.8pc/:,{\xypolygon12"A"{~={0}~>{}}},
"A1"*{\bullet},"A2"*{\bullet},"A3"*{\bullet},"A4"*{\bullet},
"A5"*{\bullet},"A6"*{\bullet},"A7"*{\bullet},"A8"*{\bullet},
"A9"*{\bullet},"A10"*{\bullet},"A11"*{\bullet},"A12"*{\bullet}\endxy
\quad \preccurlyeq \quad
\xy/r2.8pc/:,{\xypolygon12"A"{~={0}~>{}}},
"A1";"A2"**@{-},"A2";"A6"**@{-},"A6";"A1"**@{-},
"A8";"A9"**@{-},
"A1"*{\bullet},"A2"*{\bullet},"A3"*{\bullet},"A4"*{\bullet},
"A5"*{\bullet},"A6"*{\bullet},"A7"*{\bullet},"A8"*{\bullet},
"A9"*{\bullet},"A10"*{\bullet},"A11"*{\bullet},"A12"*{\bullet}\endxy
\quad \preccurlyeq \quad
\xy/r2.8pc/:,{\xypolygon12"A"{~={0}~>{}}},
"A1";"A2"**@{-},"A2";"A4"**@{-},"A4";"A6"**@{-},
"A6";"A8"**@{-},"A8";"A9"**@{-},"A9";"A1"**@{-},
"A1"*{\bullet},"A2"*{\bullet},"A3"*{\bullet},"A4"*{\bullet},
"A5"*{\bullet},"A6"*{\bullet},"A7"*{\bullet},"A8"*{\bullet},
"A9"*{\bullet},"A10"*{\bullet},"A11"*{\bullet},"A12"*{\bullet}\endxy
\quad \preccurlyeq \quad
\xy/r2.8pc/:,{\xypolygon12"A"{~={0}~>{}}},
"A1";"A2"**@{-},"A2";"A3"**@{-},"A3";"A4"**@{-},
"A4";"A5"**@{-},"A5";"A6"**@{-},"A6";"A7"**@{-},
"A7";"A8"**@{-},"A8";"A9"**@{-},"A9";"A10"**@{-},
"A10";"A11"**@{-},"A11";"A12"**@{-},"A12";"A1"**@{-},
"A1"*{\bullet},"A2"*{\bullet},"A3"*{\bullet},"A4"*{\bullet},
"A5"*{\bullet},"A6"*{\bullet},"A7"*{\bullet},"A8"*{\bullet},
"A9"*{\bullet},"A10"*{\bullet},"A11"*{\bullet},"A12"*{\bullet}\endxy
=\triv
$$
\caption{Example of a chain in $NCP(1,1,15)$}
\end{figure}

By specialising  Lemma \ref{generalncp}, we obtain:

\begin{lemma}
The poset $(NCP(1,1,n),\preccurlyeq)$ is a lattice.
\end{lemma}

Our next task is to define the \emph{complement} operation for
$NCP(1,1,n)$. This construction will later be generalised to
$NCP(e,e,n)$ (subsection \ref{subsectioncomplement}).

We start with a particular case:

\begin{definition}
\label{preccurlyeqomp}
Let $a\subset \mu_n$.
Let $x_1,\dots,x_k$ be the elements of $a$, ordered counterclockwise.
We denote by $\overline{a}$ the element of $NCP(1,1,n)$ whose parts
are
$$\rangle x_1,x_2\rangle \; , \; \rangle x_2,x_3 \rangle \; , \; \dots \; ,
\; \rangle x_{r-1},x_r \rangle \; , \; \rangle x_r,x_1\rangle.$$
\end{definition}

\begin{figure}[ht]
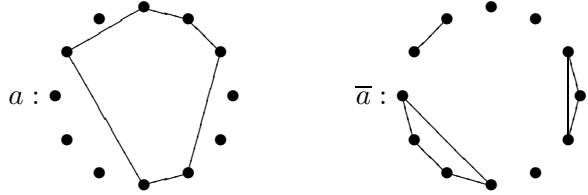

$$
a: \xy/r2.8pc/:,{\xypolygon12"A"{~={0}~>{}}},
"A2";"A3"**@{-},"A3";"A4"**@{-},"A4";"A6"**@{-},
"A6";"A10"**@{-},"A10";"A11"**@{-},"A11";"A2"**@{-},
"A1"*{\bullet},"A2"*{\bullet},"A3"*{\bullet},"A4"*{\bullet},
"A5"*{\bullet},"A6"*{\bullet},"A7"*{\bullet},"A8"*{\bullet},
"A9"*{\bullet},"A10"*{\bullet},"A11"*{\bullet},"A12"*{\bullet}\endxy
\qquad  \qquad
\overline{a}:\xy/r2.8pc/:,{\xypolygon12"A"{~={0}~>{}}},
"A5";"A6"**@{-},"A7";"A8"**@{-},"A8";"A9"**@{-},
"A9";"A10"**@{-},"A10";"A7"**@{-},
"A12";"A1"**@{-},"A1";"A2"**@{-},"A2";"A12"**@{-},
"A1"*{\bullet},"A2"*{\bullet},"A3"*{\bullet},"A4"*{\bullet},
"A5"*{\bullet},"A6"*{\bullet},"A7"*{\bullet},"A8"*{\bullet},
"A9"*{\bullet},"A10"*{\bullet},"A11"*{\bullet},"A12"*{\bullet}\endxy
$$
\caption{The map $a\mapsto \overline{a}$, when $a\subseteq \mu_n$}
\label{figureoverlinea}
\end{figure}

The general situation is as follows:
Let $u,v\in NCP(1,1,n)$, with $u\preccurlyeq v$.
Let us recursively define an element $u\backslash v\in NCP(1,1,n)$
by induction on the
number of non-singleton parts of $u$:
\begin{itemize}
\item[(i)] If $u$ has only singleton parts, we set $u\backslash v:=v$.
\item[(ii)] If $u$ has only one non-singleton part $a$, denote by $b$ the unique
part of $v$ such that $a\subseteq b$. Let $\phi$ a cyclic-order preserving
bijection from $b$ to $\mu_{|b|}$. Using Definition \ref{preccurlyeqomp}, we
obtain
an element $\overline{\phi(a)} \in NCP_{\mu_{|b|}}$. Transporting it via $\phi$,
we obtain a non-crossing partition $\phi^{-1}(\overline{\phi(a)})$ of $b$.
We set $u\backslash v$ to be the partition obtained from $v$ by splitting
the part $b$ into $\phi^{-1}(\overline{\phi(a)})$.
\item[(iii)] If $u$ has several non-singleton parts, choose $a$ one of them.
Let $u'\in NCP(1,1,n)$ be the partition obtained from $u$ by splitting
$a$ into isolated points, let $u''\in NCP(1,1,n)$ with $a$ as only 
non-singleton part. We set $u\backslash v := u' \backslash (u''\backslash v)$.
\end{itemize}

For Step (iii) to make sense, one has to observe that
$u'\preccurlyeq u''\backslash v$, which is clear (look at figure
\ref{figureoverlinea}).
The procedure is non-deterministic, since
one has to make choices when applying Step (iii). However, one may easily check
that the final result does not depend on these choices.

An obvious property is that, when $u\preccurlyeq v$, we also
have $(u\backslash v) \preccurlyeq v$.

\begin{definition}
Let $u,v\in NCP(1,1,n)$ such that $u\preccurlyeq v$.
The non-crossing partition $u\backslash v$ is called the {\em complement
of $u$ in $v$}.

Generalising Definition \ref{preccurlyeqomp}, we set,
for all $u\in NCP(1,1,n)$, $\overline{u} := u\backslash \triv$,
where $\triv$ is the maximal element of $NCP(1,1,n)$.
\end{definition}

An alternative but equivalent definition for the complement is given
in \cite{bdm}. The original construction of the complement actually
goes back to Kreweras \cite{kreweras}.

\begin{figure}[ht]
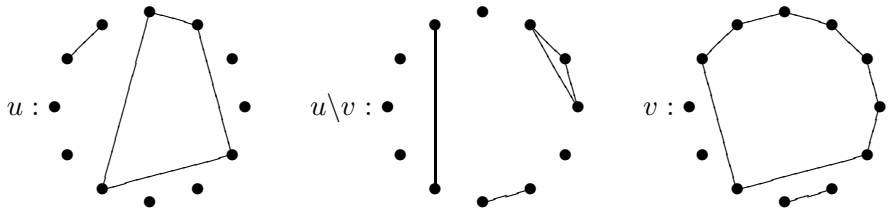

$$u:
\xy/r3pc/:,{\xypolygon12"A"{~={0}~>{}}},
"A12";"A3"**@{-},"A3";"A4"**@{-},"A4";"A9"**@{-},
"A9";"A12"**@{-},
"A5";"A6"**@{-},
"A1"*{\bullet},"A2"*{\bullet},"A3"*{\bullet},"A4"*{\bullet},
"A5"*{\bullet},"A6"*{\bullet},"A7"*{\bullet},"A8"*{\bullet},
"A9"*{\bullet},"A10"*{\bullet},"A11"*{\bullet},"A12"*{\bullet}\endxy
\qquad
u\backslash v:
\xy/r3pc/:,{\xypolygon12"A"{~={0}~>{}}},
"A1";"A2"**@{-},"A2";"A3"**@{-},"A3";"A1"**@{-},
"A9";"A5"**@{-},"A10";"A11"**@{-},
"A1"*{\bullet},"A2"*{\bullet},"A3"*{\bullet},"A4"*{\bullet},
"A5"*{\bullet},"A6"*{\bullet},"A7"*{\bullet},"A8"*{\bullet},
"A9"*{\bullet},"A10"*{\bullet},"A11"*{\bullet},"A12"*{\bullet}\endxy
\qquad
v: \xy/r3pc/:,{\xypolygon12"A"{~={0}~>{}}},
"A12";"A1"**@{-},"A1";"A2"**@{-},
"A2";"A3"**@{-},"A3";"A4"**@{-},"A4";"A5"**@{-},
"A6";"A9"**@{-},
"A9";"A12"**@{-},
"A5";"A6"**@{-},"A10";"A11"**@{-},
"A1"*{\bullet},"A2"*{\bullet},"A3"*{\bullet},"A4"*{\bullet},
"A5"*{\bullet},"A6"*{\bullet},"A7"*{\bullet},"A8"*{\bullet},
"A9"*{\bullet},"A10"*{\bullet},"A11"*{\bullet},"A12"*{\bullet}\endxy
$$
\caption{The complement: an example}
\end{figure}

\begin{figure}[ht]
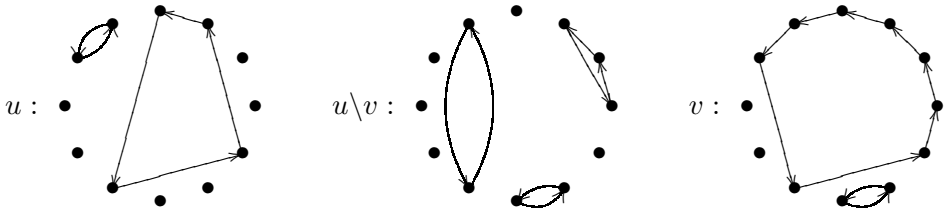

$$u:
\xy/r3pc/:,{\xypolygon12"A"{~={0}~>{}}},
{\xypolygon24"B"{~:{(1.2,0):}~>{}}},
{\xypolygon24"C"{~:{(0.8,0):}~>{}}},
"A12";"A3"**\dir{-}?>*\dir{>},"A3";"A4"**\dir{-}?>*\dir{>},
"A4";"A9"**\dir{-}?>*\dir{>},
"A9";"A12"**\dir{-}?>*\dir{>},
"A5";"A6"**\crv{"B10"}?>*\dir2{>},
"A6";"A5"**\crv{"C10"}?>*\dir2{>},
"A1"*{\bullet},"A2"*{\bullet},"A3"*{\bullet},"A4"*{\bullet},
"A5"*{\bullet},"A6"*{\bullet},"A7"*{\bullet},"A8"*{\bullet},
"A9"*{\bullet},"A10"*{\bullet},"A11"*{\bullet},"A12"*{\bullet}\endxy
\qquad
u\backslash v:
\xy/r3pc/:,{\xypolygon12"A"{~={0}~>{}}},
{\xypolygon24"B"{~:{(1.2,0):}~>{}}},
{\xypolygon24"C"{~:{(0.8,0):}~>{}}},
"A1";"A2"**\dir{-}?>*\dir{>},"A2";"A3"**\dir{-}?>*\dir{>},
"A3";"A1"**\dir{-}?>*\dir{>},
"A9";"A5"**\crv{"A0"}?>*\dir{>},
"A5";"A9"**\crv{"A7"}?>*\dir{>},
"A10";"A11"**\crv{"B20"}?>*\dir2{>},
"A11";"A10"**\crv{"C20"}?>*\dir2{>},
"A1"*{\bullet},"A2"*{\bullet},"A3"*{\bullet},"A4"*{\bullet},
"A5"*{\bullet},"A6"*{\bullet},"A7"*{\bullet},"A8"*{\bullet},
"A9"*{\bullet},"A10"*{\bullet},"A11"*{\bullet},"A12"*{\bullet}\endxy
\qquad
v: \xy/r3pc/:,{\xypolygon12"A"{~={0}~>{}}},
{\xypolygon24"B"{~:{(1.2,0):}~>{}}},
{\xypolygon24"C"{~:{(0.8,0):}~>{}}},
"A12";"A1"**\dir{-}?>*\dir{>},"A1";"A2"**\dir{-}?>*\dir{>},
"A2";"A3"**\dir{-}?>*\dir{>},"A3";"A4"**\dir{-}?>*\dir{>},
"A4";"A5"**\dir{-}?>*\dir{>},
"A6";"A9"**\dir{-}?>*\dir{>},
"A9";"A12"**\dir{-}?>*\dir{>},
"A5";"A6"**\dir{-}?>*\dir{>},
"A10";"A11"**\crv{"B20"}?>*\dir2{>},
"A11";"A10"**\crv{"C20"}?>*\dir2{>},
"A1"*{\bullet},"A2"*{\bullet},"A3"*{\bullet},"A4"*{\bullet},
"A5"*{\bullet},"A6"*{\bullet},"A7"*{\bullet},"A8"*{\bullet},
"A9"*{\bullet},"A10"*{\bullet},"A11"*{\bullet},"A12"*{\bullet}\endxy
$$
\caption{The crucial property of the complement: the composition
of the motion associated with $u$ with the motion associated
with $u\backslash v$ is homotopic to the motion associated with $v$.}
\label{figurecrucial}
\end{figure}

{\bf \flushleft Explanation.} Though there is no logical need for this yet,
interpreting the ``complement'' operation in terms of permutations and braids
provides useful intuition. To illustrate this,
we briefly recall, without proofs, how
this construction is done \cite{bdm} (this anticipates
on what will be done in Section \ref{section2} for $G(e,e,n+1)$).
To any element $u\in NCP(1,1,n)$, we associate the permutation of $\sigma_u\in
\mu_n$
sending any $\zeta\in\mu_n$ to $\zeta'$, the successor of $\zeta$ for
the counterclockwise ordering of the part of $u$ containing $\zeta$.
If we compose permutations the way we compose paths ($\sigma\tau$ meaning
``$\sigma$ then $\tau$''), we may check that the relation
$\sigma_u\sigma_{u\backslash v} = \sigma_v$, for all $u\preccurlyeq v$ in
$NCP(1,1,n)$. Since the map $u\mapsto \sigma_u$ is injective, $u\backslash v$
is uniquely determined by the equation
$\sigma_u\sigma_{u\backslash v} = \sigma_v$.
More remarkably, the analog relation holds in the usual braid
group on $n$ strings, when one associates to $u$ the braid $b_u$ represented by
motion of the points of $\mu_n$ where each $\zeta$ goes to the $\zeta'$
constructed above, following at constant speed the affine segment
$[\zeta,\zeta']$ (with the convention that, if the part of $\zeta$ and
$\zeta'$ contains only two elements, then the strings of $\zeta$
and $\zeta'$ avoid each other by ``driving on the right'' along
$[\zeta,\zeta']$). A graphical illustration is provided in Figure
\ref{figurecrucial}.

Note that the group $\mu_{n}$ naturally acts on $NCP(1,1,n)$, by its
multiplicative action on the underlying set $\mu_n$.
The map $u\mapsto \overline{u}$ is a ``square root'' of the multiplication
by $\zeta_n$:

\begin{lemma}
\label{barbar}
For all $u\in NCP(1,1,n)$, we have $\overline{\overline{u}} = \zeta_n u$.
\end{lemma}

\begin{proof}
Though it is possible to prove this in a purely combinatorial way,
we give a simple proof using the
interpretation in terms of permutations.
One has $\sigma_{u}\sigma_{\overline{u}}=\sigma_{\overline{u}}
\sigma_{\overline{\overline{u}}}=\sigma_{\triv}$,
thus $\sigma_{u}\sigma_{\triv}=\sigma_{u}\sigma_{\overline{u}}
\sigma_{\overline{\overline{u}}}=
\sigma_{\triv}\sigma_{\overline{\overline{u}}}$ and
$\sigma_{\overline{\overline{u}}}=\sigma_{\triv}^{-1}\sigma_{u}\sigma_{\triv}$.
We observe
that $\sigma_{\triv}$ is an $n$-cycle, corresponding to a rotation.
Conjugating by this $n$-cycle amounts to relabelling the underlying
set by rotation.
\end{proof}

It is well-known that the cardinality of $NCP(1,1,n)$ is
the Catalan number $\frac{1}{n+1}\binom{2n}{n}$.
We will have a natural interpretation for the cardinality of
our generalised non-crossing partitions lattices.

\begin{definition}
Let $u\in NCP(1,1,n)$. The \emph{height of $u$} is the integer
$ht(u):= n-m$, where $m$ is the number of parts of $u$.
\end{definition}

The following lemma follows straightforwardly from the definitions:
\begin{lemma}
\label{lemmahta}
For all $u,v\in NCP(1,1,n)$, we have
\begin{itemize}
\item $0\leq ht(u) \leq n-1$, $ht(\disc)=0$, $ht(\triv)=n-1$;
\item $u\preccurlyeq v \Rightarrow ht(u)\leq ht(v)$
and 
$u\prec v \Rightarrow ht(u) < ht(v)$;
\item if $u \preccurlyeq v$, then $ht(u) + ht(u\backslash v) = ht(v)$;
in particular, $ht(u)+ht(\overline{u})=n-1$.
\end{itemize}
\end{lemma}

\subsection{Non-crossing partitions of type $(e,1,n)$}
For any divisor $e$ of $n$, the group $\mu_e$ acts on $NCP(1,1,n)$ by
multiplication on the underlying set $\mu_n$.
We denote by $NCP(1,1,n)^{\mu_e}$ the set of non-crossing partitions
fixed by the action of $\mu_e$.
The following definition was already implicit in \cite{bdm}.

\begin{definition}
For any positive integers $e,n$, we set $NCP(e,1,n):=NCP(1,1,en)^{\mu_e}$.
\end{definition}

\begin{lemma}
The poset $(NCP(e,1,n),\preccurlyeq)$ is a lattice.
\end{lemma}

The reflection group $G(2,1,n)$ is also known as the ``hyperoctahedral group''
$W(B_n)$. When $e=2$, the above
definition coincides with Reiner's type $B_n$ non-crossing
partitions \cite{reiner}, and with the lattice of simple elements 
in the type $B_n$ dual braid monoid \cite{dualmonoid}.

Let $u,v\in NCP(e,1,n)$, with $u\preccurlyeq v$. The element $u\backslash v$
(defined earlier in $NCP(1,1,en)$)
is clearly in $NCP(e,1,n)$.

\begin{definition}
An element of $NCP(e,1,n)$ is {\em long} if $0$ is in the convex hull of
one of its parts (by non-crossedness, this part must be unique; it is
referred to as {\em the long part} of the partition). An element of
$NCP(e,1,n)$ which is not long is {\em short}.
\end{definition}

\begin{figure}[ht]
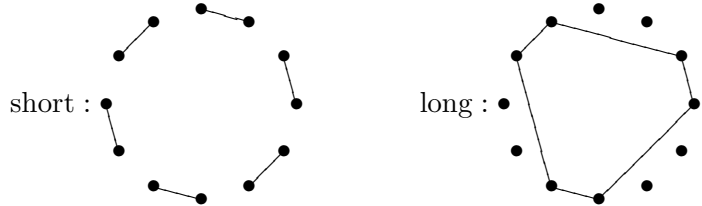

$$
\text{short}: \xy/r3pc/:,{\xypolygon12"A"{~={0}~>{}}},
"A1";"A2"**@{-},"A5";"A6"**@{-},"A9";"A10"**@{-},
"A3";"A4"**@{-},"A7";"A8"**@{-},"A11";"A12"**@{-},
"A1"*{\bullet},"A2"*{\bullet},"A3"*{\bullet},"A4"*{\bullet},
"A5"*{\bullet},"A6"*{\bullet},"A7"*{\bullet},"A8"*{\bullet},
"A9"*{\bullet},"A10"*{\bullet},"A11"*{\bullet},"A12"*{\bullet}\endxy
\qquad \qquad
\text{long}:
\xy/r3pc/:,{\xypolygon12"A"{~={0}~>{}}},
"A1";"A2"**@{-},"A2";"A5"**@{-},"A5";"A6"**@{-},
"A6";"A9"**@{-},"A9";"A10"**@{-},"A10";"A1"**@{-},
"A1"*{\bullet},"A2"*{\bullet},"A3"*{\bullet},"A4"*{\bullet},
"A5"*{\bullet},"A6"*{\bullet},"A7"*{\bullet},"A8"*{\bullet},
"A9"*{\bullet},"A10"*{\bullet},"A11"*{\bullet},"A12"*{\bullet}\endxy
$$
\caption{Some elements of $NCP(3,1,5)$}
\end{figure}

\begin{lemma}
\label{viceversa}
The map $u\mapsto \overline{u}$ sends long elements of $NCP(e,1,n)$
to short elements, and vice-versa.
In particular, the number of long elements equals the number of short
elements.

The complement of a short element in a short element is short, the
complement of a long element in a long element is short, the complement
of a short element in a long element is long.
\end{lemma}

\begin{proof}
Use the height function $ht:NCP(1,1,en)\rightarrow \BN$.
We have $ht(u) = en-m_u$, where $m_u$ is the number of parts
of $u\in NCP(e,1,n)$.

If $u$ is short, 
then $m_u\equiv 0 [e]$. If $u$ is long, we have $m_u \equiv 1 [e]$
(only the long part is fixed by $\mu_e$-action).
We have (Lemma \ref{lemmahta})
$ht(u)+ht(\overline{u})=en-1=2en-m_u-m_{\overline{u}}$. Among
$u$ and $\overline{u}$, one is short and the other is long.

The remaining statements are proved similarly.
\end{proof}

{\bf \flushleft Remark.} When studying $NCP(e,1,n)$,  
it is natural to work with a new length function $ht'$,
defined by
$ht'(u):=n-m_u/e$ when $u$ is short and $ht'(u):=n-(m_u-1)/e$
when $u$ is long. This function takes its values in $\{0,\dots,n\}$.

\subsection{Non-crossing partitions of type $(e,e,n+1)$}

In this subsection, we define non-crossing partitions of
type $(e,e,n+1)$ as being non-crossing partitions
of $\mu_{en}\cup\{0\}$ satisfying certain additional conditions.
The geometry of $NCP(e,e,n+1)$ will be related to that of
$NCP(e,1,n)$ by three natural poset morphisms, according to the diagram:

$$\xymatrix{ NCP(e,1,n) \ar@{^{(}->}[rr]^{*} &  
& NCP(e,e,n+1) \ar@{>>}@/_1.5em/[ll]_{\sharp}
\ar@{>>}@/^1.5em/[ll]^{\flat} }$$

\begin{definition}
For all $u \in NCP_{\mu_{en}\cup\{0\}}$, 
we denote by $u^{\flat}\in NCP_{\mu_{en}}$ the partition obtained by
forgetting $0$.
We set $$NCP(e,e,n+1) := \{u\in NCP_{\mu_{en}\cup\{0\}} | u^{\flat}\in
NCP(e,1,n)\}.$$

An element $u\in NCP(e,e,n+1)$ is said to be:
\begin{itemize}
\item \emph{short symmetric},
if $u^{\flat}$ is a short element of $NCP(e,1,n)$ and
$\{0\}$ is a part in $u$;
\item \emph{long symmetric}, if $u^{\flat}$
is a long element of $NCP(e,1,n)$. The non-crossing condition
then implies that $a\cup \{0\}$ is a part of $u$, where $a$ is the
long part of $u^{\flat}$; we say that $a\cup \{0\}$ is \emph{the
long part} of $u$;
\item \emph{asymmetric}, if $u^{\flat}$ is a short
element of $NCP(e,1,n)$, and
$\{0\}$ is a not a part of $u$. There is then an unique part $a$
of $u^{\flat}$ such that $a\cup\{0\}$ is a part of $u$; the part
$a\cup\{0\}$ is \emph{the asymmetric part of $u$}.
\end{itemize}
\end{definition}

Clearly,
these three cases are mutually exclusive and any element of $NCP(e,e,n+1)$
is of one of the three types.

\begin{figure}[ht]
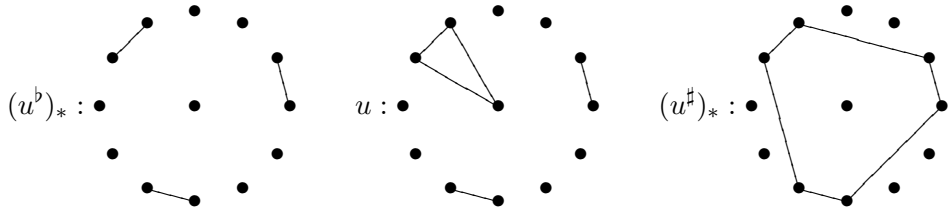

$$(u^{\flat})_*:
 \xy/r3pc/:,{\xypolygon12"A"{~={0}~>{}}},
"A1";"A2"**@{-},"A5";"A6"**@{-},"A9";"A10"**@{-},
"A0"*{\bullet},"A1"*{\bullet},"A2"*{\bullet},"A3"*{\bullet},"A4"*{\bullet},
"A5"*{\bullet},"A6"*{\bullet},"A7"*{\bullet},"A8"*{\bullet},
"A9"*{\bullet},"A10"*{\bullet},"A11"*{\bullet},"A12"*{\bullet}\endxy
\qquad
u:
\xy/r3pc/:,{\xypolygon12"A"{~={0}~>{}}},
"A1";"A2"**@{-},"A5";"A6"**@{-},"A9";"A10"**@{-},
"A0";"A6"**@{-},"A0";"A5"**@{-},
"A0"*{\bullet},"A1"*{\bullet},"A2"*{\bullet},"A3"*{\bullet},"A4"*{\bullet},
"A5"*{\bullet},"A6"*{\bullet},"A7"*{\bullet},"A8"*{\bullet},
"A9"*{\bullet},"A10"*{\bullet},"A11"*{\bullet},"A12"*{\bullet}\endxy
\qquad
(u^{\sharp})_*:
\xy/r3pc/:,{\xypolygon12"A"{~={0}~>{}}},
"A1";"A2"**@{-},"A5";"A6"**@{-},"A9";"A10"**@{-},
"A2";"A5"**@{-},"A6";"A9"**@{-},"A1";"A10"**@{-},
"A0"*{\bullet},"A1"*{\bullet},"A2"*{\bullet},"A3"*{\bullet},"A4"*{\bullet},
"A5"*{\bullet},"A6"*{\bullet},"A7"*{\bullet},"A8"*{\bullet},
"A9"*{\bullet},"A10"*{\bullet},"A11"*{\bullet},"A12"*{\bullet}\endxy
$$
\caption{Some elements of 
$NCP(3,3,5)$:
short symmetric, asymmetric and long symmetric.}
\end{figure}

From the definition, it follows immediately that if $u$ is short
symmetric, $v$ asymmetric and $w$ long symmetric, then one cannot
have $v\preccurlyeq u$ nor $w\preccurlyeq v$ nor $w\preccurlyeq u$.

The subposet of $NCP(e,e,n+1)$ consisting of 
symmetric elements, both short and long, 
is isomorphic to $NCP(e,1,n)$ via the map $NCP(e,1,n)\rightarrow NCP(e,e,n+1),
u \mapsto u_*$, with $u_*$ defined as follows:
\begin{itemize}
\item If $u$ is short, we set $u_*:=u\cup\{\{0\}\}$.
This identifies short elements in $NCP(e,1,n)$ with short symmetric
elements
in $NCP(e,e,n+1)$.
\item If $u$ is long and if $a$ is the long part of $u$, the
partition $u_*$ is obtained from $u$ by replacing $a$ by
$a\cup\{0\}$. This identifies long elements in $NCP(e,1,n)$ with
long symmetric elements in $NCP(e,e,n+1)$.
\end{itemize}

Finally, we define a map $NCP(e,e,n+1) \rightarrow NCP(e,1,n), \lambda \mapsto
\lambda^{\sharp}$ as follows:
\begin{itemize}
\item If $u$ is symmetric, then we set
$u^{\sharp}:=u^{\flat}$.
\item If $u$ is asymmetric, let $a$ be the part containing $0$.
Let $\tilde{a}:= (\bigcup_{\zeta\in\mu_e} \zeta a)-\{0\}$. We
set $u^{\sharp}$ to be the element of $NCP(e,1,n)$ containing
$\tilde{a}$ as a part, and apart from that made of parts in $u$.
\end{itemize}

From the definitions, it immediately follows that, for all $u\in NCP(e,1,n)$ and
for all $v\in NCP(e,e,n+1)$:
$$(u_*)^\flat =u \quad \text{and} \quad (u_*)^\sharp =u$$
$$(v^\flat)_* \preccurlyeq v \preccurlyeq (v^\sharp)_*$$
$$(v^\flat)_* =v \Leftrightarrow (v^\sharp)_*=v \Leftrightarrow
\text{$v$ is symmetric}.$$
Not only are the maps $\flat$ and $\sharp$ retractions of the natural inclusion
$*$, but they can also be viewed as ``adjoints'' of $*$, in the following
sense (the lemma is an easy consequence of the above formulae):

\begin{lemma}
\label{adjoint}
For all $u \in NCP(e,1,n)$, for all $v\in NCP(e,e,n+1)$, we
have
$$v \preccurlyeq u_* \Leftrightarrow v^{\sharp} \preccurlyeq u \qquad \text{and}
\qquad u_* \preccurlyeq v \Leftrightarrow u \preccurlyeq v^{\flat}.$$
\end{lemma}

As expected, we have:

\begin{lemma}
\label{latticeproperty}
The poset $(NCP(e,e,n+1),\preccurlyeq)$ is a lattice.
\end{lemma}

\begin{proof}
Once again, it is enough to check that $NCP(e,e,n+1)$ is stable by 
set-theoretical meet. Let $v,v'\in NCP(e,e,n+1)$. Let $v\wedge v'$ be their
meet in $NCP_{\mu_{en}\cup \{0\}}$.
We observe that
$(v\wedge v')^{\flat}=v^{\flat}\wedge v'^{\flat}$, which
is an element of $NCP(e,1,n)$. Thus $v\wedge v'\in NCP(e,e,n+1)$.
\end{proof}

\subsection{Type $(1,1,n+1)$ intervals in $NCP(e,e,n+1)$}
\begin{definition}
For $\zeta\in \mu_{en}$, we set
$$a_{\zeta}:=\{0,\zeta_{en}^1\zeta,\zeta_{en}^2\zeta,
\dots,\zeta_{en}^n\zeta\}.$$
and we consider the following asymmetric elements
of $NCP(e,e,n+1)$:
\begin{itemize}
\item the element $m_\zeta$, the partition with the unique non-singleton part
$\{0,\zeta\}$;
\item the element $M_\zeta$, made up of asymmetric part
$a_{\zeta}$,
the remaining parts being of the form
$\zeta'a_{\zeta}-\{0\}$,
for $\zeta'\in \mu_e-\{1\}$.
\end{itemize}
\end{definition}

When $u,v\in NCP(e,e,n+1)$, we use the notation
$$[u,v]:=\{w\in NCP(e,e,n+1) | u\preccurlyeq w \preccurlyeq v\}.$$
The following lemma characterises certain intervals in $NCP(e,e,n+1)$.
We denote by $\res_{\zeta}$ the restriction morphism from
the lattice of partitions of $\mu_{en}\cup \{0\}$ to the lattice
of partitions of $a_{\zeta}$. Since $a_{\zeta}$ is a strictly convex
$(n+1)$-gon, we have $NCP_{a_{\zeta}}\simeq NCP(1,1,n+1)$.
\begin{lemma}
\label{intervals}
\begin{itemize}
\item[(i)] 
The set of minimal (resp. maximal) asymmetric elements in 
$NCP(e,e,n+1)$ is $\{m_\zeta | \zeta\in\mu_{en}\}$ (resp. 
$\{M_\zeta | \zeta\in\mu_{en}\}$).
\item[(ii)] Let $\zeta\in \mu_{en}$. 
The map $\res_{\zeta}$ induces a poset isomorphism
$$\varphi_\zeta:([\disc,M_{\zeta}],\preccurlyeq)
\stackrel{\sim}{\rightarrow} (NCP_{a_{\zeta}},\preccurlyeq).$$
\item[(iii)] Let $\zeta\in\mu_{en}$.
The map $\res_{\zeta}$ induces a poset isomorphism
$$\psi_\zeta:([m_{\zeta},\triv],\preccurlyeq)
\stackrel{\sim}{\rightarrow} (NCP_{a_{\zeta}},\preccurlyeq).$$
\item[(iv)]
The composition $\psi^{-1}_{\zeta} \varphi_{\zeta}$ is
the map $u\mapsto m_{\zeta}\vee u$.
\end{itemize}
\end{lemma}

Since $a_{\zeta}$ is convex and has cardinal
$n+1$, we have $(NCP_{a_{\zeta}},\preccurlyeq)\simeq
(NCP(1,1,n+1),\preccurlyeq)$. It should also be noted that
$m_{\zeta}$ is not finer than $M_{\zeta}$, thus 
$[\disc,M_{\zeta}] \cap [m_{\zeta},\triv]= \varnothing$.

\begin{proof}
(i) is clear.

(ii) Any part $b$ of a given $u\in [\disc,M_\zeta]$ is either
finer than $a_{\zeta}$ or finer than some $\zeta'a_\zeta-\{0\}$,
$\zeta'\in\mu_e-\{1\}$. The symmetry conditions imply that
given any $v\in NCP_{a_{\zeta}}$, there is a unique way to complete
it with partitions of the $\zeta'a_\zeta-\{0\}$ to obtain an element
of $[\disc,M_\zeta]$ (explicitly, the partition of $\zeta'a_\zeta-\{0\}$ 
is the restriction of $\zeta'v$).

(iii) 
Let $b_{\zeta}:=a_{\zeta}\cup \{\zeta\}$.
Elements of $[m_\zeta,\triv]$ are either asymmetric or long symmetric.
Using the symmetry conditions,
it is clear that they are uniquely determined by their restriction
to $b_{\zeta}$, and that
 $0$ and $\zeta$ lie in the same part of this restriction.
Conversely, any element of $NCP_{b_{\zeta}}$
satisfying the condition  that $0$ and $\zeta$ are connected
may be obtained this way.
It is also clear that the restriction from $b_{\zeta}$ to $a_{\zeta}$
identifies $\{v\in NCP_{b_{\zeta}} | \text{$0$ and $\zeta$ are in the
same part} \}$ with $NCP_{a_{\zeta}}$. This proves (iii).

It is clear that $u\mapsto m_{\zeta} \vee u$ induces a map
$\theta:[\disc,M_\zeta]\rightarrow [m_{\zeta},\triv]$.
To obtain (iv), it is enough to check that, for all $u\in[\disc,M_\zeta]$,
$\res_{\zeta}(u)=\res_{\zeta}\theta(u)$, which is straightforward.
\end{proof}

\subsection{The complement}
\label{subsectioncomplement}

\begin{definition}
\label{defcomp}
Let $u,v\in NCP(e,e,n+1)$, with $u\preccurlyeq v$. We define an
element $u\backslash v$, the \emph{complement of $u$ in $v$},
as follows:
\begin{itemize}
\item[(A)] If $u$ and $v$ are both symmetric, we set
$$u\backslash v := (u^\flat \backslash v^\flat)_*.$$
\item[(B)] If $v$ is asymmetric, we choose $\zeta\in\mu_{en}$ such
that $v\preccurlyeq M_{\zeta}$.
We set
$$u\backslash v:= \varphi_{\zeta}^{-1} 
(\varphi_\zeta(u)\backslash \varphi_{\zeta}(v)),$$
where the complement operation is defined on
$NCP_{a_{\zeta}}$ via a standard identification with $NCP(1,1,n+1)$.
\item[(C)] If $u$ is asymmetric, we choose
$\zeta\in\mu_{en}$ such
that $m_{\zeta}\preccurlyeq u$.
We set
$$u\backslash v:= \varphi_{\zeta}^{-1}
(\psi_\zeta(u)\backslash \psi_{\zeta}(v)),$$
where the complement operation is defined on
$NCP_{a_{\zeta}}$ via a standard identification with $NCP(1,1,n+1)$.
(Note that we really mean $\varphi^{-1}$ and not $\psi^{-1}$. This implies that
 $u\backslash v$ lies in $[\disc,M_\zeta]$).
\end{itemize}
\end{definition}

For the above definition to make sense, one has to check that,
in case (B) and (C), the result does not depend on the choice of $\zeta$;
also, when both $u$ and $v$ are asymmetric, one has to check that (B)
and (C) give the same $u\backslash v$. We leave these to the reader
to verify.

As it was the case in $NCP(1,1,n)$, this definition becomes more
natural when interpreted in the associated braid group
(see Proposition \ref{mappingb} below).

An essential property of the complement is that
for all $u\preccurlyeq v$, one has
$$u\backslash v\preccurlyeq v.$$
Indeed, for the cases (A) and (B) of the definition, this follows immediately
from the corresponding results for classical non-crossing partitions.
In the situation (C), we have $\psi_\zeta(u)\backslash \psi_{\zeta}(v)
\preccurlyeq \psi_{\zeta}(v)$. Since $\varphi_{\zeta}$ is a poset morphism,
we have $u\backslash v \preccurlyeq
\varphi^{-1}_{\zeta}\psi_{\zeta}(v)$. Looking at the construction
of $\varphi_\zeta$, one observes that, for any
$w\in NCP_{a_{\zeta}}$, for any $\tilde{w}\in NCP(e,e,n+1)$ such
that $\res_{a_{\zeta}}(\tilde{w})=w$, one has $\varphi^{-1}_\zeta(w)
\preccurlyeq \tilde{w}$. This applies in particular to $\tilde{w}=v$.
Thus $u\backslash v \preccurlyeq
\varphi^{-1}_{\zeta}\psi_{\zeta}(v) \preccurlyeq v$.

We also use the following notation:

\begin{definition}
For all $u\in NCP(e,e,n+1)$, we set $\overline{u}:= u \backslash \triv$.
\end{definition}

For example, since $m_{\zeta}$ is the minimal element of $[m_{\zeta},\triv]$,
$\psi_\zeta(m_{\zeta})$ is the discrete element of $NCP_{a_{\zeta}}$; clearly,
$\psi_\zeta(\triv)$ is the trivial element $\{a_{\zeta}\}$. Thus
$\psi_\zeta(m_{\zeta}) \backslash \psi_\zeta(\triv)
=\psi_\zeta(m_{\zeta}) \backslash \{a_{\zeta}\} = \{a_{\zeta}\}$, which
maps to $M_{\zeta}$ by $\varphi_\zeta^{-1}$.
Thus $$\overline{m_{\zeta}}=M_{\zeta}.$$
To compute $M_{\zeta}$, there are several possible $\zeta'$ such
that $m_{\zeta'}\preccurlyeq M_{\zeta}$. Take for example
$\zeta':=\zeta_{en} \zeta$.
Then $\psi_{\zeta'}(M_\zeta)$ is the partition of $a_{\zeta'}$ with 
two parts, $\{0,\zeta_{en}^2\zeta,\dots,\zeta_{en}^n\zeta\}$ and
$\{\zeta_{en}^{n+1}\zeta\}$. Thus $\psi_{\zeta'}(M_\zeta)\backslash
\psi_{\zeta'}(\triv)= \psi_{\zeta'}(M_\zeta)\backslash
\{a_{\zeta'}\}$ is the partition with $\{0,\zeta_{en}^{n+1}\zeta\}$
as only non-singleton part. The latter partition is mapped
to $m_{\zeta_{en}^{n+1}\zeta}= \zeta_{en}^{n+1}
m_{\zeta}$ by $\varphi_{\zeta'}$.
We have proved that
$$\overline{M_{\zeta}}=\zeta_{en}^{n+1}m_{\zeta}.$$
As a consequence, for any $m$ (resp. $M$) minimal (resp. maximal)
asymmetric,
$$\overline{\overline{m}} = \zeta_{en}^{n+1} m \qquad
\text{and} \qquad 
\overline{\overline{M}} = \zeta_{en}^{n+1} M.$$

We list below some basic properties of the complement:
\begin{lemma}
\label{lemmacon}
\begin{itemize}
\item[(i)] The map $u\mapsto \overline{u}$ is a poset anti-automorphism
of $NCP(e,e,n+1)$.
\item[(ii)] If $u$ is short symmetric, then $\overline{u}$ is long symmetric.
\item[(iii)] If $u$ is long symmetric, then $\overline{u}$ is short symmetric.
\item[(iv)] If $u$ is asymmetric, then $\overline{u}$ is asymmetric.
\item[(v)] For all $u$, $\overline{\overline{u}}=\zeta_e \zeta_{en} u$.
\end{itemize}
\end{lemma}

{\bf \flushleft Remark.}
By (ii) and (iii), the map $\phi:u\mapsto \overline{\overline{u}}$ is
an automorphism of $NCP(e,e,n+1)$ whose order is the order of
$\zeta_e\zeta_{en}= \zeta_{en}^{n+1}$, that is, $\frac{en}{(n+1)\wedge e}$.
It is worth noting that there are several types of orbits for this
automorphism:
\begin{itemize}
\item By definition, symmetric partitions are preserved under multiplication
by $\mu_e$, so the action of $\phi$ on symmetric partitions is by
multiplication by $\zeta_{en}$; this action has order $n$ (the action
is not free; e.g., $\disc$ and $\triv$ are fixed points).
\item The group $\mu_{en}$ acts freely on the set of asymmetric partitions,
which decomposes, under the action of $\phi$, into orbits of equal cardinal
$\frac{en}{(n+1)\wedge e}$.
\end{itemize}

\subsection{Height}
Let $u\in NCP(e,e,n+1)$. Let $m$ be the number of parts of $u$.
If $u$ is symmetric, then parts of $u$ come in orbits of cardinal
$e$ for the action of $\mu_e$, except the isolated part $\{0\}$ (when
$u$ is short) and the long part (when $u$ is long). In both symmetric cases,
$m \equiv 1 [e]$. When $u$ is asymmetric, $e|m$.

We define the \emph{height} $ht(u)$ of $u$
as follows:
\begin{itemize}
\item If $u$ is short symmetric, we
set $ht(u) := n - (m-1)/e$.
\item If $u$ is asymmetric, we set
$ht(u) := n+1 - m/e$.
\item If $u$ is long symmetric, we
set $ht(u) := n+1 - (m-1)/e$.
\end{itemize}

It is not difficult to check the basic properties listed below:
\begin{lemma}
\label{height}
Let $u\in NCP(e,e,n+1)$.
\begin{itemize}
\item[(i)] If $u$ is short symmetric, then $0\leq ht(u) \leq n-1$.
\item[(ii)] If $u$ is asymmetric, then $1\leq ht(u) \leq n$.
\item[(iii)] If $u$ is long symmetric, then $2\leq ht(u) \leq n+1$.
\item[(iv)] For all $v\in NCP(e,e,n+1)$ such that $u\preccurlyeq v$, we have
$ht(u)+ht(u\backslash v) = ht(v)$. In particular,
$ht(u)+ht(\overline{u}) = n+1$.
\item[(v)] For all $v\in NCP(e,e,n+1)$, we have $u\preccurlyeq v 
\Rightarrow ht(u) \leq ht(v)$ and $u\prec v \Rightarrow ht(u) < ht(v)$.
\item[(vi)] For all $v\in NCP(e,e,n+1)$ such that $u\preccurlyeq
v$, we may find a chain
$$u=u_0\preccurlyeq u_1 \preccurlyeq \dots \preccurlyeq u_k=v$$ in
$NCP(e,e,n+1)$ such that, for all $i$ in $0,\dots,k$, we have
$ht(u_i)= ht(u)+i$.
\end{itemize}
\end{lemma}

\section{Simple elements}
\label{section2}

\subsection{The complex reflection group $G(e,e,n+1)$}
View $\FS_{n+1}$ in its natural representation as the group of
$(n+1)\times (n+1)$ permutation matrices.
For any positive integer $d$,
let $\Delta(de,e,n+1)\subset
\GL_{n+1}(\BC)$ be the group of diagonal matrices with diagonal
coefficients in $\mu_{de}$ and determinant in $\mu_d$.
The complex reflection group $G(de,e,n+1)$ is, by definition,
the subgroup of $\GL_{n+1}(\BC)$ generated by $\FS_{n+1}$ and
$\Delta(de,e,n+1)$ (see for example \cite{bmr} for a general 
description of the classification of complex reflection groups).
Then
$$G(de,e,n+1) \simeq \Delta(de,e,n+1) \rtimes \FS_{n+1}.$$

The group $G(e,e,n+1)$ has $e^n (n+1)!$ elements, among them are
$e(n+1)(n+2)/2$ reflections, all with order $2$ and forming
a single conjugacy class. A minimal set of generating reflections
has cardinal $n+1$: take for example the $n$ permutation matrices
$\sigma_i:=(i \; \;  i+1)$ (for $i=1,\dots,n$) plus the reflection
$\Diag(\zeta_e,\overline{\zeta_e},1,1,\dots, 1,1) \sigma_1$.

Let $X_0,\dots,X_n$ be the canonical coordinates on $V:=\BC^{n+1}$
(we consider that matrices start with a $0$-th line and a $0$-th column).
The reflecting hyperplanes have equations $X_i = \zeta X_j$, where 
$i,j\in\{0,\dots,n\}$ and $\zeta\in \mu_e$. They form a single orbit under the
action of $G(e,e,n+1)$.

It is sometimes useful to compare $G(e,e,n+1)$ with other groups in the
infinite series of complex reflection groups.
By definition, we have
$$G(e,e,n+1) \subset G(e,1,n+1).$$
We may also construct a monomorphism
$$G(e,1,n) \stackrel{\psi}{\hookrightarrow} G(e,e,n+1)$$ 
as follows. Let $\chi$ the character of $G(e,1,n)$ trivial
on $\FS_n$ and coinciding with the determinant on $\Delta(e,1,n)$. Then,
for any matrix $M\in G(e,1,r)$, the matrix (given by blocks)
$$\psi(M):=\left( \begin{matrix} \overline{\chi(g)} & 0 \\ 0 & M
\end{matrix} \right) $$
is in $G(e,e,n+1)$.

\subsection{The braid group $B(e,e,n+1)$}
\label{Been}

Let $V:=\BC^{n+1}$, endowed with the action of $G(e,e,n+1)$.
Let $V^{\reg}$ be the complement in $V$ of the union of the reflection
hyperplanes. By definition, the braid group $B(e,e,n+1)$ associated
with $G(e,e,n+1)$ is the fundamental group of the quotient space
$V^{\reg}/G(e,e,n+1)$. This definition involves the choice of a
basepoint. Our preferred basepoint is chosen by means of Springer
theory of regular elements (\cite{springer}).

The invariant degrees are
$e,2e,\dots,ne,n+1$. The largest degree, $ne$, is regular, in the sense
of Springer (this means that
there exists a element of $G(e,e,n+1)$ of order $ne$ and having an
eigenvector in $V^{\reg}$; such an element is said to be \emph{regular}).
A typical regular element of order $ne$ is
$$c:=
\begin{pmatrix}
\overline{\zeta_e}  & 0 & 0 & \dots & \dots & 0 & 0 \\
0 & 0 & 1 & 0 & \dots & \dots & 0    \\
0 & 0 & 0 & 1 &  0 & \dots & 0   \\
\vdots  & \vdots & & \ddots & \ddots & & \vdots \\
\vdots  & \vdots &  & & \ddots & \ddots & \vdots\\
0 & 0 & 0 & \dots  & 0 & 0 & 1 \\
0 & \zeta_e & 0 & \dots & \dots & 0 & 0
\end{pmatrix}
$$
An example of regular eigenvector in $\ker(c-\zeta_{ne}\Id)$ is
$$x_0:=\begin{pmatrix} 0 \\ \zeta_{ne} \\ \zeta_{ne}^2 \\ \vdots \\
\zeta_{ne}^n=\zeta_e\end{pmatrix}.$$

We choose $x_0$ as basepoint for $V^{\reg}$ and its orbit $\overline{x_0}$
as basepoint for the quotient space $V^{\reg}/G(e,e,n+1)$.
Hence we fix the (now fully explicit) definition:
$$B(e,e,n+1):=\pi_1(V^{\reg}/G(e,e,n+1), \overline{x_0}).$$

{\bf \flushleft Remarks.}
\begin{itemize}
\item
The initial symmetry between all reflecting
hyperplanes is broken by the choice of $c$. This explains why our
presentation involves two distinct types of generators.
\item
Note that $c\in \psi(G(e,1,n))$. The element $c$ is actually the
image via $\psi$ of a maximal regular element in $G(e,1,n)$.
A consequence of the main results in \cite{bdm}
is that we have a Garside monoid for $B(e,1,n)$, related to the choice
of a maximal regular element, and admitting a description using
planar partitions of $\mu_{en}$ invariant by $\mu_e$.
\end{itemize}

\subsection{From non-crossing partitions to elements of $B(e,e,n+1)$
and $G(e,e,n+1)$}
\label{defmappingb}

We construct here a map $NCP(e,e,n+1) \rightarrow B(e,e,n+1)$.

We work with a given $u\in NCP(e,e,n+1)$.
To simplify notations, we sometimes 
omit referring to $u$, though the construction
of course depends on $u$.

For every $z\in \mu_{ne}\cup\{0\}$, denote by $z'$
the successor of $z$ for the counterclockwise ordering of the part of $u$
containing $z$ (when $u$ is long symmetric, when dealing with the long part,
which is not convex, we consider the counterclockwise ordering of the 
non-zero elements of this part, and set $0':=0$).

We consider a path $\gamma_{u,z}$ (also denoted by $\gamma_z$ if there
is no ambiguity on $u$) as follows:
\begin{itemize}
\item If the part of $u$ containing $z$ contains three or more elements,
we set
\begin{eqnarray*}
\gamma_z: [0,1] & \longrightarrow & \BC \\
t & \longmapsto & (1-t) z + t z'
\end{eqnarray*}
\item If the part of $u$ containing $z$ contains exactly two elements
($z$ and $z'$), we  slightly perturbate the above
definition by setting
\begin{eqnarray*}
\gamma_z: [0,1] & \longrightarrow & \BC \\
t & \longmapsto & (1-t) z + t z'+it(1-t)(z-z') \varepsilon_n
\end{eqnarray*}
where $i$ is the standard square root of $-1$ used to fix the orientation of
$\BC$, and $\varepsilon_n\in\BR_{>0}$ is a fixed ``small enough'' number,
depending only on $n$.
\item If $z$ is isolated in $u$,
we take for $\gamma_z$ the constant path $t\mapsto z$.
\end{itemize}

\begin{figure}[ht]
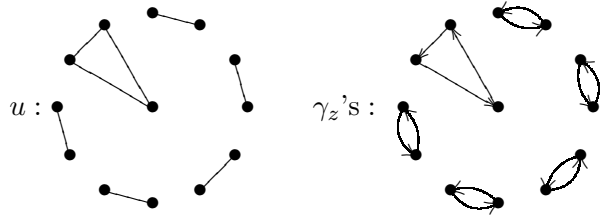

$$u:
\xy/r3pc/:,{\xypolygon12"A"{~={0}~>{}}},
"A1";"A2"**@{-},"A5";"A6"**@{-},"A9";"A10"**@{-},
"A0";"A6"**@{-},"A0";"A5"**@{-},
"A3";"A4"**@{-},"A7";"A8"**@{-},"A11";"A12"**@{-},
"A0"*{\bullet},"A1"*{\bullet},"A2"*{\bullet},"A3"*{\bullet},"A4"*{\bullet},
"A5"*{\bullet},"A6"*{\bullet},"A7"*{\bullet},"A8"*{\bullet},
"A9"*{\bullet},"A10"*{\bullet},"A11"*{\bullet},"A12"*{\bullet}\endxy
\qquad \text{$\gamma_z$'s}:
\xy/r3pc/:,{\xypolygon12"A"{~={0}~>{}}},
{\xypolygon24"B"{~:{(1.2,0):}~>{}}},
{\xypolygon24"C"{~:{(0.8,0):}~>{}}},
"A0";"A5"**\dir{-}?>*\dir{>},"A5";"A6"**\dir{-}?>*\dir{>},
"A6";"A0"**\dir{-}?>*\dir{>},
"A11";"A12"**\crv{"B22"}?>*\dir2{>},
"A12";"A11"**\crv{"C22"}?>*\dir2{>},
"A3";"A4"**\crv{"B6"}?>*\dir2{>},
"A4";"A3"**\crv{"C6"}?>*\dir2{>},
"A1";"A2"**\crv{"B2"}?>*\dir2{>},
"A2";"A1"**\crv{"C2"}?>*\dir2{>},
"A9";"A10"**\crv{"B18"}?>*\dir2{>},
"A10";"A9"**\crv{"C18"}?>*\dir2{>},
"A7";"A8"**\crv{"B14"}?>*\dir2{>},
"A8";"A7"**\crv{"C14"}?>*\dir2{>},
"A0"*{\bullet},"A6"*{\bullet},"A5"*{\bullet},"A11"*{\bullet},
"A12"*{\bullet},"A1"*{\bullet},"A2"*{\bullet},"A3"*{\bullet},
"A4"*{\bullet},"A7"*{\bullet},"A8"*{\bullet},"A9"*{\bullet},"A10"*{\bullet}
\endxy
$$
\caption{The various $\gamma_z$'s, for a given $u\in NCP(3,3,5)$.}
\label{examplezu}
\end{figure}

We use these paths in $\BC$ to associate to $u$ a path in $V^{\reg}$,
with initial point $x_0$. When $u$ is asymmetric,
there is unique $z\in \mu_{ne}$ such that $z'=0$ (in the example
illustrated above, $z=\zeta_{12}^5$).
We say that $j\in \{1,\dots, n\}$ is \emph{special} (with respect to $u$)
if $u$ is asymmetric and $j$ is such that
$z\in \mu_e \zeta_{ne}^j$, for the unique $z$ such that $z'=0$
(when $u$ is asymmetric, this condition selects a unique $j\in \{1,\dots,n\}$;
when $u$ is symmetric, there are no special integers; in our above example,
the special integer is $1$).

\begin{figure}[ht]
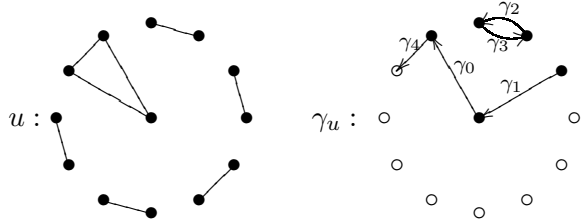

$$u:
\xy/r3pc/:,{\xypolygon12"A"{~={0}~>{}}},
"A1";"A2"**@{-},"A5";"A6"**@{-},"A9";"A10"**@{-},
"A0";"A6"**@{-},"A0";"A5"**@{-},
"A3";"A4"**@{-},"A7";"A8"**@{-},"A11";"A12"**@{-},
"A0"*{\bullet},"A1"*{\bullet},"A2"*{\bullet},"A3"*{\bullet},"A4"*{\bullet},
"A5"*{\bullet},"A6"*{\bullet},"A7"*{\bullet},"A8"*{\bullet},
"A9"*{\bullet},"A10"*{\bullet},"A11"*{\bullet},"A12"*{\bullet}\endxy
\qquad \gamma_u:
\xy/r3pc/:,{\xypolygon12"A"{~={0}~>{}}},
{\xypolygon24"B"{~:{(1.2,0):}~>{}}},
{\xypolygon24"C"{~:{(0.8,0):}~>{}}},
"A0";"A5"**\dir{-}?>*\dir{>} ?(.4)*!LD!/-5pt/^{_{\gamma_0}},
"A5";"A6"**\dir{-}?>*\dir{>} ?(.6)*!LD!^{^{\gamma_4}},
"A2";"A0"**\dir{-}?>*\dir{>}  ?(.6)*!LD!^{^{\gamma_1}},
"A3";"A4"**\crv{"B6"}?>*\dir2{>} ?(.4)*!LD!/-5pt/^{_{\gamma_2}},
"A4";"A3"**\crv{"C6"}?>*\dir2{>},
"C6"*{_{\gamma_3}},
"A0"*{\bullet},"A2"*{\bullet},"A3"*{\bullet},"A4"*{\bullet},
"A5"*{\bullet},"A6"*{\circ},"A7"*{\circ},"A8"*{\circ},
"A9"*{\circ},"A10"*{\circ},"A11"*{\circ},"A12"*{\circ},"A1"*{\circ}\endxy
$$
\caption{Illustration of the map $u\mapsto \gamma_u$; the initial
positions of the $\gamma_j$'s are marked by black dots.}
\label{examplebu}
\end{figure}

If $j\in \{1,\dots,n\}$ is special, we set $\gamma_j:=\zeta^{-1} \gamma_{z}$,
where $z$ is the above predecessor of $0$, and $\zeta$ is the element
of $\mu_e$ such that $z=\zeta \zeta_{ne}^j$.
If $j$ is not special, we set $\gamma_j:=\gamma_{\zeta_{ne}^j}.$
The path $\gamma_u$ is defined by
\begin{eqnarray*}
\gamma_u: [0,1] & \longrightarrow & \BC^n \\
t & \longmapsto & \begin{pmatrix} \gamma_{0}(t) \\ \gamma_1(t) \\
\vdots \\ \gamma_{n}(t)\end{pmatrix} 
\end{eqnarray*}
Figure \ref{examplebu} illustrates this construction.

Provided that $\varepsilon_n$ is taken close enough to $0$, this
indeed defines a path in $V^{\reg}$, whose homotopy class does not
depend on the explicit $\varepsilon_n$.

The final endpoint $\gamma_u(1)$
of $\gamma_u$ lies in the orbit $\overline{x_0}$.
This follows from the following easy lemma, left to the reader:

\begin{lemma}
\label{orbit}
A vector $(z_0,\dots,z_n)\in \BC^{n+1}$  is in the orbit
$\overline{x_0}$ if and only if we have
$\{z_0,\dots,z_n\} = X \cup \{0\}$, where $X\subset \mu_{ne}$ is a
set of coset representatives of $\mu_{ne}/\mu_e$.
\end{lemma}

Thus the path $\gamma_u$ defines a loop $\overline{\gamma_u}$
in the quotient $V^{\reg}/G(e,e,n+1)$.

\begin{definition}
We denote by $b_u$ the element of $B(e,e,n+1)$
represented by the loop $\overline{\gamma_u}$.
We set $$P_B:=\{b_u | u \in NCP(e,e,n+1)\}.$$
\end{definition}

As in \cite{bmr}, we have a fibration exact sequence
$$\xymatrix@1{1\ar[r]&P(e,e,n+1)\ar[r] & B(e,e,n+1)\ar[r]^{\pi}
& G(e,e,n+1)\ar[r] & 1},$$
where $P(e,e,n+1)=\pi_1(V^{\reg},x_0)$.
The interpretation of the morphism $\pi$ is as follows: given $b\in B(e,e,n+1)$,
we may represent it by a loop in $V^{\reg}/G(e,e,n+1)$ with endpoint
$\overline{x_0}$; this loops lifts to a unique path in $V^{\reg}$ with
initial point $x_0$; the final endpoint $x$ of this path is in the orbit
$\overline{x_0}$; the element $\pi(b)$ is the unique $g\in G(e,e,n+1)$
such that $g(x_0)=x$.

\begin{definition}
For all $u\in NCP(e,e,n+1)$, we set $g_u:=\pi(b_u)$.
We set $$P_G:=\pi(P_B) = \{g_u | u \in NCP(e,e,n+1)\}.$$
\end{definition}

So far we have constructed a commutative diagram
$$\xymatrix{
 & P_B \ar[dd]^{\pi}
 \\
NCP(e,e,n+1) \ar[ur]^{u\mapsto b_u} \ar[dr]_{u\mapsto g_u} \\
& P_G  }$$
Our task in the next sections will be to prove that these maps are
bijective -- actually, $P_G$ will be given
a more intrinsic definition and a natural poset structure, making
$u\mapsto g_u$ a poset isomorphism.

We may endow $P_B$ with the \emph{left divisibility partial order},
again denoted by $\preccurlyeq$ and defined by
$$\forall b,b' \in P_B, b\preccurlyeq b' \stackrel{\text{def}}
{\Longleftrightarrow}
\exists b''\in P_B, bb''=b'.$$
In this setting, a consequence of the next proposition 
is that the map $u\mapsto b_u$ is a poset morphism (it will later be proved
to be a poset isomorphism).

\begin{proposition}
\label{mappingb}
For any $u,v\in NCP(e,e,n+1)$ such that $u\preccurlyeq v$, we have
$b_u b_{u\backslash v} = b_v$.
\end{proposition}

\begin{figure}[ht]
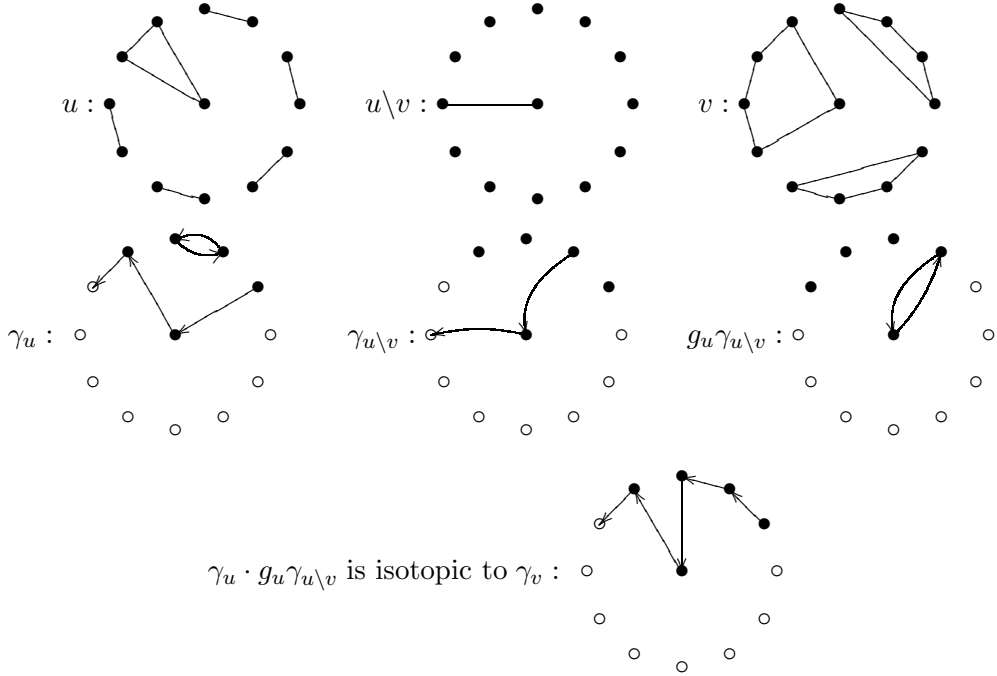

$$u:
\xy/r3pc/:,{\xypolygon12"A"{~={0}~>{}}},
"A1";"A2"**@{-},"A5";"A6"**@{-},"A9";"A10"**@{-},
"A0";"A6"**@{-},"A0";"A5"**@{-},
"A3";"A4"**@{-},"A7";"A8"**@{-},"A11";"A12"**@{-},
"A0"*{\bullet},"A1"*{\bullet},"A2"*{\bullet},"A3"*{\bullet},"A4"*{\bullet},
"A5"*{\bullet},"A6"*{\bullet},"A7"*{\bullet},"A8"*{\bullet},
"A9"*{\bullet},"A10"*{\bullet},"A11"*{\bullet},"A12"*{\bullet}\endxy
\qquad 
u\backslash v:
\xy/r3pc/:,{\xypolygon12"A"{~={0}~>{}}},
"A0";"A7"**@{-},
"A0"*{\bullet},"A1"*{\bullet},"A2"*{\bullet},"A3"*{\bullet},"A4"*{\bullet},
"A5"*{\bullet},"A6"*{\bullet},"A7"*{\bullet},"A8"*{\bullet},
"A9"*{\bullet},"A10"*{\bullet},"A11"*{\bullet},"A12"*{\bullet}\endxy
\qquad
v:
\xy/r3pc/:,{\xypolygon12"A"{~={0}~>{}}},
"A1";"A2"**@{-},"A5";"A6"**@{-},"A9";"A10"**@{-},"A6";"A7"**@{-},
"A0";"A8"**@{-},"A0";"A5"**@{-},
"A2";"A3"**@{-},"A1";"A4"**@{-},
"A3";"A4"**@{-},"A7";"A8"**@{-},"A11";"A12"**@{-},
"A10";"A11"**@{-},"A9";"A12"**@{-},
"A0"*{\bullet},"A1"*{\bullet},"A2"*{\bullet},"A3"*{\bullet},"A4"*{\bullet},
"A5"*{\bullet},"A6"*{\bullet},"A7"*{\bullet},"A8"*{\bullet},
"A9"*{\bullet},"A10"*{\bullet},"A11"*{\bullet},"A12"*{\bullet}\endxy
$$
$$
\gamma_u:
\xy/r3pc/:,{\xypolygon12"A"{~={0}~>{}}},
{\xypolygon24"B"{~:{(1.2,0):}~>{}}},
{\xypolygon24"C"{~:{(0.8,0):}~>{}}},
"A0";"A5"**\dir{-}?>*\dir{>},
"A5";"A6"**\dir{-}?>*\dir{>}, 
"A2";"A0"**\dir{-}?>*\dir{>},  
"A3";"A4"**\crv{"B6"}?>*\dir2{>}, 
"A4";"A3"**\crv{"C6"}?>*\dir2{>},
"A0"*{\bullet},"A2"*{\bullet},"A3"*{\bullet},"A4"*{\bullet},
"A5"*{\bullet},"A6"*{\circ},"A7"*{\circ},"A8"*{\circ},
"A9"*{\circ},"A10"*{\circ},"A11"*{\circ},"A12"*{\circ},"A1"*{\circ}\endxy
\qquad
\gamma_{u\backslash v}:
\xy/r3pc/:,{\xypolygon12"A"{~={0}~>{}}},
{\xypolygon12"C"{~:{(0.5,0):}~>{}}},
"A0";"A7"**\crv{"C6"}?>*\dir{>},
"A3";"A0"**\crv{"C4"}?>*\dir{>},
"A0"*{\bullet},"A2"*{\bullet},"A3"*{\bullet},"A4"*{\bullet},
"A5"*{\bullet},"A6"*{\circ},"A7"*{\circ},"A8"*{\circ},
"A9"*{\circ},"A10"*{\circ},"A11"*{\circ},"A12"*{\circ},"A1"*{\circ}\endxy
\qquad
g_u\gamma_{u\backslash v}:
\xy/r3pc/:,{\xypolygon12"A"{~={0}~>{}}},
{\xypolygon12"C"{~:{(0.5,0):}~>{}}},
"A0";"A3"**\crv{"C2"}?>*\dir{>},
"A3";"A0"**\crv{"C4"}?>*\dir{>},
"A0"*{\bullet},"A2"*{\circ},"A3"*{\bullet},"A4"*{\bullet},
"A5"*{\bullet},"A6"*{\bullet},"A7"*{\circ},"A8"*{\circ},
"A9"*{\circ},"A10"*{\circ},"A11"*{\circ},"A12"*{\circ},"A1"*{\circ}\endxy
$$
$$
\gamma_u \cdot g_u\gamma_{u\backslash v}\; \text{is isotopic to}\; \gamma_v:
\xy/r3pc/:,{\xypolygon12"A"{~={0}~>{}}},
{\xypolygon24"B"{~:{(1.2,0):}~>{}}},
{\xypolygon24"C"{~:{(0.8,0):}~>{}}},
{\xypolygon12"D"{~:{(0.5,0):}~>{}}},
"A0";"A5"**\dir{-}?>*\dir{>},
"A5";"A6"**\dir{-}?>*\dir{>}, 
"A2";"A3"**\dir{-}?>*\dir{>},  
"A3";"A4"**\dir{-}?>*\dir{>}, 
"A4";"A0"**\dir{-}?>*\dir{>}, 
"A0"*{\bullet},"A2"*{\bullet},"A3"*{\bullet},"A4"*{\bullet},
"A5"*{\bullet},"A6"*{\circ},"A7"*{\circ},"A8"*{\circ},
"A9"*{\circ},"A10"*{\circ},"A11"*{\circ},"A12"*{\circ},"A1"*{\circ}\endxy
$$
\caption{Illustration of the relation $b_ub_{u\backslash v}=b_v$ on an
example.
The black dots indicate the coordinates
at the starting point of the paths.
The product $b_ub_{u\backslash v}$ is represented in $V^{\reg}$
by the path obtained by concatenating $\gamma_u$ and $g_u \cdot
\gamma_{u\backslash v}$; this concatenation is homotopic
to $\gamma_v$.} 
\label{examplepropbu}
\end{figure}

Note that we have already seen a result of this nature, for the type $A$
situation (see Figure \ref{figurecrucial}).
Figure \ref{examplepropbu} illustrates
the relation with an example in $B(3,3,5)$. The paths $\gamma_u$ and
$\gamma_{u\backslash v}$ may not be concatenated in $V^{\reg}$; one
insteads concatenates $\gamma_u$ (whose final point is $g_u(x_0)$)
with the transformed path $g_u \gamma_{u\backslash v}$.

\begin{proof}
We have to deal with the three cases of Definition \ref{defcomp}.

(A) In \cite{bdm}, $NCP(e,1,n)$ was interpreted in relation with
the braid group $B(e,1,n)$ of $G(e,1,n)$ (the latter being viewed as
the centraliser of the $e$-th root of the full twist in $G(1,1,en)$).
It is easy to see that our definition of $u\backslash v$ and $b_u$,
for symmetric elements of $NCP(e,e,n+1)$, coincide via $\flat$ with
the corresponding constructions from \emph{loc. cit.}. The desired result
is then a consequence of \emph{loc. cit.}, Proposition 1.8 (ii) (the
element there denoted by $\delta_{\lambda}$ is our $b_u$).
Another way to recover the case (A), without using \cite{bdm}, is by a direct
check using Sergiescu relations.

(B) Assume $v$ is asymmetric, and choose $\zeta$ such that
$v\preccurlyeq M_\zeta$.
Instead of representing $b_v$ by the path $\gamma_v$, which starts at $x_0$,
we may replace $\gamma_u$ by any $g(\gamma_u)$ with $g\in G(e,e,n+1)$.
Thanks to Lemma \ref{orbit}, we may choose $g$ such that
$$g(x_0) = \begin{pmatrix} 0 \\ \zeta_{en}^1 \zeta \\ \vdots
\\ \zeta_{en}^n\zeta \end{pmatrix}.$$
A direct computation shows that
$$g(\gamma_v) = \begin{pmatrix} \gamma_{v,0} \\ \gamma_{v,\zeta_{en}^1 \zeta}
\\ \vdots
\\ \gamma_{v,\zeta_{en}^n\zeta} \end{pmatrix}.$$
Since $v\preccurlyeq M_\zeta$, we have
$$\{\gamma_{v,0}(0),\gamma_{v,\zeta_{en}^1 \zeta}(0),\dots,
\gamma_{v,\zeta_{en}^n\zeta}(0)\} = 
\{\gamma_{v,0}(1),\gamma_{v,\zeta_{en}^1 \zeta}(1),\dots,
\gamma_{v,\zeta_{en}^n\zeta}(1)\} = a_{\zeta}.$$ In particular, the
final endpoint $g(\gamma_v)(1)$ is related to the starting point
$g(\gamma_v)(0)$ by a permutation matrix $\sigma\in G(1,1,n+1)$.
In other words, $g(\gamma_v)$ may be viewed as a classical type $A$ braid
involving $n+1$ strings, with initial and final position at $a_{\zeta}$.
A similar discussion applies to $u$. In this classical type $A$ braid group,
there is a complement operation and the relation similar to the one we
want to prove here (see Figure \ref{figurecrucial}, or \cite{bdm}).
Definition \ref{defcomp} is precisely designed to be compatible with
this classical relation: the classical relation asserts that two braids
are homotopic; one may easily sees that the homotopy may be chosen so that
at all time, all strings lie in the convex hull of $a_{\zeta}$; this
implies that the desired homotopy holds in $V^{\reg}$.

We leave the details to the reader, as well as the proof in case (C), which
may be carried out in a similar manner.
\end{proof}

\section{Non-crossing partitions of height $1$ and local generators}
\label{section3}

In this section,
we have a closer look at the image by the previous maps
of the set of non-crossing
partitions of height $1$. The image by $u\mapsto b_u$ will be the generating set
in our new presentation for $B(e,e,n+1)$.

\begin{definition}
\label{h1defi}
Let $p,q\in \BZ$ with $0 < |p-q| < n$. We
set $u_{p,q}$ to be the (short symmetric) element of $NPC(e,e,n+1)$
whose non-singleton parts are
$$\{\zeta_{ne}^p,\zeta_{ne}^q\}\; ,\;
\zeta_e\{\zeta_{ne}^p,\zeta_{ne}^q\}\; ,\; \zeta_e^2
\{\zeta_{ne}^p,\zeta_{ne}^q\}\; , \; \dots \; 
,\; \zeta_e^{e-1}\{\zeta_{ne}^p,\zeta_{ne}^q\}.$$
We extend this notation to the situation where $p$ and $q$
are the images in $\BZ/en\BZ$ of integers $\tilde{p}$ and $\tilde{q}$
with $0 < |\tilde{p}-\tilde{q}| < n$.

For all $p\in \BZ$, we set
$u_p$ to be the (asymmetric) element of $NCP(e,e,n+1)$ with
only non-singleton part
$\{0,\zeta_{en}^{p}\}$. We extend this notation to the situation
where $p\in \BZ/en\BZ$.
\end{definition}

\begin{figure}[ht]
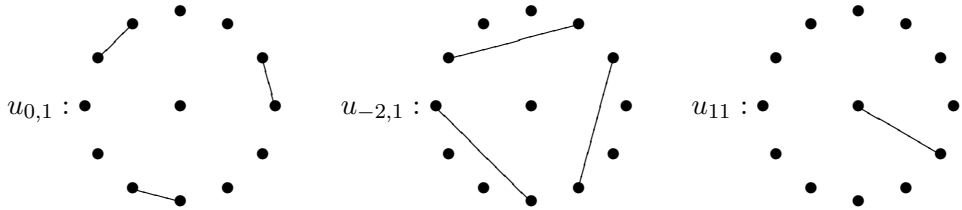

$$
u_{0,1}: \xy/r3pc/:,{\xypolygon12"A"{~={0}~>{}}},
"A1";"A2"**@{-},"A5";"A6"**@{-},"A9";"A10"**@{-},
"A0"*{\bullet},"A1"*{\bullet},"A2"*{\bullet},"A3"*{\bullet},"A4"*{\bullet},
"A5"*{\bullet},"A6"*{\bullet},"A7"*{\bullet},"A8"*{\bullet},
"A9"*{\bullet},"A10"*{\bullet},"A11"*{\bullet},"A12"*{\bullet}\endxy
\qquad
u_{-2,1}:
\xy/r3pc/:,{\xypolygon12"A"{~={0}~>{}}},
"A3";"A6"**@{-},"A7";"A10"**@{-},"A11";"A2"**@{-},
"A0"*{\bullet},"A1"*{\bullet},"A2"*{\bullet},"A3"*{\bullet},"A4"*{\bullet},
"A5"*{\bullet},"A6"*{\bullet},"A7"*{\bullet},"A8"*{\bullet},
"A9"*{\bullet},"A10"*{\bullet},"A11"*{\bullet},"A12"*{\bullet}\endxy
\qquad
u_{11}:
\xy/r3pc/:,{\xypolygon12"A"{~={0}~>{}}},
"A0";"A12"**@{-},
"A0"*{\bullet},"A1"*{\bullet},"A2"*{\bullet},"A3"*{\bullet},"A4"*{\bullet},
"A5"*{\bullet},"A6"*{\bullet},"A7"*{\bullet},"A8"*{\bullet},
"A9"*{\bullet},"A10"*{\bullet},"A11"*{\bullet},"A12"*{\bullet}\endxy
$$
\caption{Some height $1$ elements of $NCP(3,3,5)$}
\end{figure}

\begin{lemma}
\label{listh1}
There are $en(n-1)$ unordered pairs of elements $p,q\in\BZ/en\BZ$ who
have representatives $\tilde{p},\tilde{q}\in \BZ$ with $0<|\tilde{p}-
\tilde{q} | < n$. For all such $p,q$, we have
$u_{p,q}=u_{p+n,q+n}$. Conversely, if $u_{p',q'}=u_{p,q}$,
then $p'=p+kn,q'=q+kn$ or $p'=q+kn,q'=p+kn$ for some $k$.
Thus there are $n(n-1)$ distinct partitions of type $u_{p,q}$.

There are $en$ distinct partitions of type $u_p$.

The $n(n-1)$ elements of the form $u_{p,q}$, together with the
$en$ elements of the form $u_{p}^{(i)}$, is a complete list
of height $1$ elements of $NCP(e,e,n+1)$.
\end{lemma}

In view of Lemma \ref{height} (v) and (vi), minimal height $1$ non-crossing
partitions are precisely minimal non-discrete elements in $NCP(e,e,n+1)$.

\begin{definition}
With $p,q$ corresponding to the situations considered
in Definition \ref{h1defi}, we set
$$a_{p,q}:=b_{u_{p,q}} \qquad \text{and} \qquad
a_{p}:=b_{u_p}.$$

We set $A:=\{a_{p,q} | p,q\in \BZ, 0<|p- q|<n \} \cup
\{a_{p}| p\in \BZ/en\BZ\}$.
\end{definition}

Since all reflections in $G(e,e,n+1)$ have order $2$, the
notions of {\em visible hyperplanes}, {\em local generators}
and {\em local monoids} may be imported at no cost from section
3 of \cite{dualmonoid} to our context.
The following proposition states that the geometric interpretation
of the dual braid monoid carries on to our situation.
It is not difficult to prove the following result -- this is 
just a matter of basic affine
geometry -- but since we do not use it later we do not include a proof.

\begin{proposition}
The elements of $A$ are 
braid reflections (aka ``meridiens'' or ``generators-of-the-monodromy'').
The corresponding reflecting hyperplanes are visible
from our basepoint $x_0$, and $A$ is actually the set
of local generators at $x_0$.
No other hyperplanes are visible from $x_0$. Hence $A$
generate the local braid monoid at $x_0$.
\end{proposition}

For finite real reflection groups, among local generators, one may always
find
``classical generators'' (corresponding to an Artin presentation);
in particular, local generators indeed generate the braid group
(\cite{dualmonoid}, Propositions 3.4.3(1) and 3.4.5).
In our situation, one may wonder whether a Brou\'e-Malle-Rouquier generating
set may be found inside $A$.
This may be done; however first, 
for the convenience of the reader, let us quote from
\cite{bmr} their result about $B(e,e,n+1)$:

\begin{theorem}[Brou\'e-Malle-Rouquier]
\label{bmr1}
The group $B(e,e,n+1)$ admits the following presentation:
\noindent
\emph{Generators:} $\tau_2, \tau'_2, \tau_3, \ldots, \tau_n, \tau_{n+1}$

\noindent
\emph{Relations:}  The commuting relations are: 
$\tau_i \tau_j = \tau_j \tau_i$ whenever $|i-j| \geq 2$,
together with $\tau'_2 \tau_j = \tau_j \tau'_2$ for all $j \geq 4$.
The others are:
$$\begin{array}{rcll}
\tau_i \tau_j &= &\tau_j \tau_i &\mbox{ for } |i-j| \geq 2\\
\tau'_2 \tau_j &= &\tau_j \tau'_2 &\mbox{ for }j \geq 4\\
\langle \tau_2 \tau'_2 \rangle^e &= &\langle \tau'_2 \tau_2 \rangle^e \\
\tau_i \tau_{i+1} \tau_i &= &\tau_{i+1} \tau_i \tau_{i+1} &\mbox{ for }i = 2, \ldots, n\\
\tau'_2 \tau_3 \tau'_2 &= &\tau_3 \tau'_2 \tau_3\\
\tau_3 \tau_2 \tau'_2 \tau_3 \tau_2 \tau'_2 &= &\tau_2 \tau'_2 \tau_3  \tau_2 \tau'_2 \tau_3
\end{array}$$
where the expression $\langle ab \rangle^k$ denotes the alternating product
$aba\cdots$ of length $k$.
\end{theorem}

Brou\'e-Malle-Rouquier construct their generators (which are braid reflections)
in an inductive manner,
using suitable fibrations. The actual generators depend on a certain
number of choices, and by having a careful look at the proof in \cite{bmr},
it is not difficult to rephrase the above theorem in the more precise manner:

\begin{theorem}[after Brou\'e-Malle-Rouquier]
\label{bmr2}
Denote by $BMR$ the abstract group presented by the above Brou\'e-Malle-Rouquier
presentation. There is an isomorphism
\begin{eqnarray*}
BMR & \stackrel{\sim}{\longrightarrow} & B(e,e,n+1) \\
\tau_2 & \longmapsto & a_0 \\
\tau'_2 & \longmapsto & a_n \\
\tau_i & \longmapsto & a_{i-3,i-2} \quad  \text{for $i=3,\dots,n+1$}. 
\end{eqnarray*}
\end{theorem}

\begin{figure}[ht]
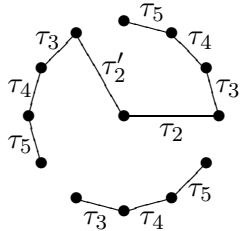

$$
\xy/r3pc/:,{\xypolygon12"A"{~={0}~>{}}},
"A1";"A2"**@{-} ?*_!/-6pt/{\tau_3},"A5";"A6"**@{-} ?*_!/-6pt/{\tau_3},
"A9";"A10"**@{-} ?*_!/-6pt/{\tau_3},
"A3";"A2"**@{-} ?*_!/6pt/{\tau_4},"A7";"A6"**@{-} ?*_!/6pt/{\tau_4},
"A11";"A10"**@{-} ?*_!/6pt/{\tau_4},
"A3";"A4"**@{-} ?*_!/-6pt/{\tau_5},"A7";"A8"**@{-} ?*_!/-6pt/{\tau_5},
"A11";"A12"**@{-} ?*_!/-6pt/{\tau_5},
"A0";"A1"**@{-} ?*_!/-6pt/{\tau_2},"A0";"A5"**@{-} ?*_!/-6pt/{\tau'_2},
"A0"*{\bullet},"A1"*{\bullet},"A2"*{\bullet},"A3"*{\bullet},"A4"*{\bullet},
"A5"*{\bullet},"A6"*{\bullet},"A7"*{\bullet},"A8"*{\bullet},
"A9"*{\bullet},"A10"*{\bullet},"A11"*{\bullet},"A12"*{\bullet}\endxy
$$
\caption{Typical Brou\'e-Malle-Rouquier generators for $B(3,3,5)$}
\end{figure}

We complete this section with a computation of the reflecting hyperplanes
of the non-crossing reflections. In view of Lemma \ref{listh1},
it is not difficult to see that the
the considered indices fully parametrise the height $1$ non-crossing
partitions.

\begin{lemma}
\label{hyper1}
Let $p,q,\in \BZ$.

If $1\leq p < q \leq n$, then $g_{u_{p,q}}$
is the reflection with hyperplane $X_p=X_q$, and
$g_{u_{q,p+n}}$ is the reflection with hyperplane $\zeta_eX_p=X_q$.

If $1 \leq p \leq en$, write $p=in+j$, with $j\in\{1,\dots,n\}$. Then
$g_{u_{p}}$ is the reflection with hyperplane $X_0=\zeta_e^iX_j$.
\end{lemma}

\begin{proof}
The computation goes as explained in Subsection \ref{mappingb}, from where
we take our notation.

If $1\leq p < q \leq n$,  the final point of $\gamma_{u_{p,q}}$ is
$$(0,\zeta_{en},\dots,\zeta_{en}^{p-1},\zeta_{en}^{q},
\zeta_{en}^{p+1},\dots,\zeta_{en}^{q-1},\zeta_{en}^{p},
\zeta_{en}^{q+1},\dots,\zeta_{en}^n).$$ Thus $g_{u_{p,q}}$ is the permutation
matrix associated with the transposition $(p \; q)$; in other words, it
is the reflection with hyperplane $X_p=X_q$.
The final point of $\gamma_{u_{q,p+n}}$ is
$$(0,\zeta_{en},\dots,\zeta_{en}^{p-1},\zeta_e^{-1}\zeta_{en}^{q},
\zeta_{en}^{p+1},\dots,\zeta_{en}^{q-1},\zeta_e\zeta_{en}^{p},
\zeta_{en}^{q+1},\dots,\zeta_{en}^n).$$ The element of $G(e,e,n+1)$ sending
$x_0$ to this final point resembles the permutation matrix of $(p \; q)$,
except that the submatrix of lines $p,q$ and columns $p,q$ is
$$\begin{pmatrix}
0 & \zeta_e^{-1} \\
\zeta_e & 0
\end{pmatrix}$$
This is the reflection with hyperplane $\zeta_eX_p=X_q$.

If $i\in \{0,\dots,e-1\}$, $j\in \{1,\dots,n\}$, the
final point of $\gamma_{u_{in+j}}$ is
$$(\zeta_{en}^{in+j},\zeta_{en},\dots,\zeta_{en}^{j-1},0,
\zeta_{en}^{j+1},\dots,\zeta_{en}^n).$$
The associated matrix resembles the permutation matrix of $(0 \; p)$,
except that the submatrix of lines $0,j$ and columns $0,j$ is
$$\begin{pmatrix}
0 & \zeta_e^i \\
\zeta_{e}^{-i} & 0
\end{pmatrix}$$
(to understand where the $\zeta_{e}^{-i}$ comes from, remember that
the product of the non-zero coefficients must be $1$).
This is the reflection with hyperplane $X_0=\zeta_e^iX_j$.
\end{proof}

\section{A length function on $G(e,e,n+1)$}
\label{section4}

Let $T$ be the set of all reflections in $G(e,e,n+1)$. As in Section 1 of
\cite{dualmonoid}, we consider, for each $g\in G(e,e,n+1)$, the set
$\Red_T(g)$ of \emph{reduced $T$-decompositions} of $g$, \ie, of
minimal length sequences $(t_1,\dots,t_k)$ such that $g=t_1\dots t_k$, and
we denote by $l_T(g)$ the length of such minimal sequences.
We also have a partial order relation $\preccurlyeq_T$: for all
$g,h\in G(e,e,n+1)$, $g\preccurlyeq_T h$ if and only if there exists 
$(t_1,\dots,t_k)\in \Red_T(h)$ and $l\leq k$ such
that $g=t_1\dots t_l$.
(Since $T$ is a union of conjugacy classes, one obtains an equivalent
definition when replacing the condition $g=t_1\dots t_l$ by
$g=t_{i_1}\dots t_{i_l}$, for some increasing sequence
$1\leq i_1< i_2 < \dots < i_l \leq k$.)

For obvious reasons from elementary linear algebra,
we have a first approximation for the function $l_T$:
$$\forall g\in G(e,e,n+1), l_T(g) \geq \codim(\ker(g-1)).$$

A major difference with the situation with real reflection groups is
that the above inequality may indeed be strict, as we will see shortly.
The analog of the 
basic Lemma 1.3.1 of \cite{dualmonoid} does not hold  for all $g$
but, as it will
appear shortly, many results from {\em loc. cit.} continue to hold.

In the following lemma, $c$ is the element defined in 
Subsection \ref{Been}.

\begin{lemma}
\label{lemmadivisors}
\begin{itemize}
\item[(i)] We have $c=g_{\triv}$.
\item[(ii)]
For all $u\in NCP(e,e,n+1)$, we have $ht(u)=\codim(\ker(g_u-1))=l_T(g_u)$.
\item[(iii)]
For all $u\in NCP(e,e,n+1)$, we have
$$\ker(g_u-1) = \bigcap_{v \preccurlyeq u, ht(v)=1} \ker(g_{v}-1).$$
\end{itemize}
\end{lemma}

\begin{proof} Part (i) is obtained directly from the definitions.

Note that the case $ht(u)=1$ of (ii) follows from the fact
that $g_u$ is then a reflection.

Let $u\in NCP(e,e,n+1)$.
We prove first the inequality $ht(u) \geq l_T(u)$ from
(ii). By Lemma \ref{height}(vi), we may find a ``filtration''
$u_0\preccurlyeq u_1
\preccurlyeq u_2 \preccurlyeq \dots \preccurlyeq u_{ht(u)}$ with
$u_{ht(u)}=u$ and $\forall i, ht(u_i)=i$. For $i=1,\dots,ht(u)$,
set $v_i:=u_{i-1}\backslash u_i$. By Lemma \ref{height}(iv), we
have $ht(v_i)=1$ and, as we have seen above, $g_{v_i}$ is a reflection.
Since
$$g_u=g_{u_{ht(u)}}= g_{u_{ht(u)-1}} g_{v_{ht(u)}}=
g_{u_{ht(u)-2}}g_{v_{ht(u)-1}}g_{v_{ht(u)}}
= \dots =g_{v_1}g_{v_2}\dots g_{v_{ht(u)}},$$
we have a $T$-decomposition of length $ht(u)$ of $u$.

The case $u=\triv$ of (ii) is a direct check.

We may now prove (ii) in full generality.
For any $u\in NCP(e,e,n+1)$, we have $c=g_ug_{\overline{u}}$
(by Proposition \ref{mappingb}) and
$ht(c)=ht(u)+ht(\overline{u})$ (by Lemma \ref{height}(iv)).
Thus
$$ht(u)+ht(\overline{u})=ht(c)=
\codim(\ker(c-1)) = \codim(\ker(g_ug_{\overline{u}}-1)).$$
By basic linear algebra, $\ker(g_ug_{\overline{u}}-1)
\supseteq \ker(g_u-1) \cap \ker(g_{\overline{u}}-1)$
and $\codim(\ker(g_u-1) \cap \ker(g_{\overline{u}}-1))
\leq \codim(\ker(g_u-1)) + \codim(\ker(g_{\overline{u}}-1))$.
Thus
$$ht(u)+ht(\overline{u}) \leq
\codim(\ker(g_u-1)) + \codim(\ker(g_{\overline{u}}-1)).$$
But on the other hand, we already checked that
$$ht(u) \geq l_T(g_u) \geq \codim(\ker(g_u-1))$$ and
$$ht(\overline{u}) \geq l_T(g_{\overline{u}}) \geq
\codim(\ker(g_{\overline{u}}-1)).$$
All considered inequalities must be equalities, and so (ii) follows.

Note that, given as earlier a ``filtration''
$u_0\preccurlyeq u_1
\preccurlyeq u_2 \preccurlyeq \dots \preccurlyeq u_{ht(u)}$ with
$u_{ht(u)}=u$ and $\forall i, ht(u_i)=i$, the same basic linear
algebra considerations, together with (ii), imply that
$\ker(g_u-1) = \bigcap_{i=1}^n \ker(g_{v_i}-1)$ (where $v_i:=u_{i-1}^{-1}u_i$).
Since, by Lemma
\ref{height}(vi), we may start this filtration with any partition $u_1$
of height $1$
such that $u_1\preccurlyeq u$, 
we obtain (iii).
\end{proof}

\begin{proposition}
\label{precTdiv}
Let $u,v\in NCP(e,e,n+1)$. 
The following assertions are equivalent:
\begin{itemize}
\item[(i)] $u\preccurlyeq v$,
\item[(ii)] $g_u\preccurlyeq_T g_v$,
\item[(iii)] $\ker (g_u-1) \supseteq \ker(g_v-1)$.
\end{itemize}
\end{proposition}

\begin{proof}
Assuming (i), one has $g_u g_{u\backslash v} =g_v$, with
$l_T(g_u) + l_T(g_{u\backslash v}) =  ht(u) + ht(u\backslash v) =
ht(v) = l_T(g_v)$ (using Lemma \ref{lemmadivisors}
(ii) and Lemma \ref{height} (iv)), thus (ii).

Assuming (ii), we have $(t_1,\dots,t_k)\Red_T(g_v)$ and $l\leq k$ such 
that $(t_1,\dots,t_l)\in\Red_T(g_u)$.
We have $\cap_{i=1}^k \ker(t_i-1) \subseteq \ker(g_u-1)$.
A priori, $\codim(\cap_{i=1}^k \ker(t_i-1)) \leq k$.
By Lemma \ref{divisors} (ii), we know that $\codim(\ker(g_u-1))=k$.
Thus $\cap_{i=1}^k \ker(t_i-1) = \ker(g_u-1)$.
Similarly, $\cap_{i=1}^l \ker(t_i-1) = \ker(g_v-1)$. This gives (iii).

To prove $\text{(iii)} \Rightarrow \text{(i)}$,
we first observe that one may use Lemma \ref{hyper1}
and Lemma \ref{lemmadivisors} to find an explicit description
of $\ker(g_v-1)$, as follows.

Let $$B_v:=\{ \bigcup_{\zeta\in\mu_e} \zeta a | 
\text{$a$ part of $v$} \}.$$
Each element of $B_v$ is a union of parts of $v$.
Clearly, $B_v$ is a (possibly crossing) partition of $\mu_{en}\cup \{0\}$,
with $v\preccurlyeq B_v$.
To any part $a$ of $B_v$, we assign an integer $i_a$ as follows.
Since $a$ is closed under $\mu_{e}$-action, 
it intersects $\{0,\zeta_{en},\zeta_{en}^2,\dots, \zeta_{en}^n\}$
(which is a fundamental domain for the $\mu_{e}$-action on $\mu_{en}\cup
\{0\}$). If $0\in a$, set $i_a:=0$. If $0\notin a$, choose $i_a\in\{1,\dots,n\}$
such that $\zeta_{en}^{i_a}\in a$.
If $v$ is short symmetric or asymmetric, we set
$$I_v:=\{ i_a | \text{$a$ part of $B_v$}\};$$
if $v$ is long symmetric, we set
$$I_v:=\{ i_a | \text{$a$ part of $B_v$}, 0\notin a\}.$$
\begin{quote}
Claim. \emph{We have $|I_v| = \dim (\ker(g_v -1))$.
For any function $\alpha:I_v\rightarrow \BC$, there exists
a unique $x\in \ker(g_v -1)$ such that for all $i\in I_v$, the $i$-th
coordinate of $x$ is $\alpha(i)$.}
\end{quote}
The first statement of the claim
follows easily from a case-by-case analysis of the 
$\mu_e$ action on the set of parts of $v$, from the defining formula
for $ht(v)$, and from the relation $ht(v)=\codim(\ker(g_v-1))$ (Lemma
\ref{lemmadivisors} (ii)).

To prove the second statement, we observe that for any
$i,j\in\{1,\dots,n\}$ such that $\zeta_{en}^i$ and $\zeta_{en}^j$ lie
in the same part of $B_v$, then $\zeta_{en}^i$ and $\zeta\zeta_{en}^j$
lie in the same part of $v$ for some $\zeta\in \mu_e$, and 
$\ker(g_v-1)$ is contained in the corresponding
hyperplane, of equation $X_i=\zeta'X_j$
for some $\zeta'\in\mu_e$ (use Lemma \ref{lemmadivisors} (iii) and
Lemma \ref{hyper1}). For the part $a_0$ of $B_v$ containing $0$,
we similarly remark that $\ker(g_v-1)$ is included in hyperplanes
of equations $X_0=\zeta''X_i$, for each $i\in\{1,\dots,n\}$ such that
$\zeta_{en}^i\in a_0$. When $v$ is long symmetric, we note in addition
that the long part contains both $0,\zeta_{en}^i,\zeta_{en}^{i+n}$ for
some $i$, thus we have (once again by Lemma \ref{lemmadivisors} (iii) and
Lemma \ref{hyper1}) relations $X_0=X_i=\zeta_e X_i$, and $X_0=0$.
Altogether, these observations imply that any $x\in \ker(g_v-1)$ is entirely
determined by its $I_v$ coordinates. 
One concludes using the relation $|I_v| = \dim(\ker(g_v-1))$.

We are now ready to prove $\text{(iii)}\Rightarrow \text{(i)}$. Assume that 
(i) does not hold. Choose a generic $x$ in $\ker(g_v-1)$.
Since $u$ is not finer than $v$, we may find
$\zeta,\zeta'\in \mu_{en}\cup \{0\}$ lying in the same part of $u$, but not
in the same part of $v$. Up to multiplying by an element of $\mu_{e}$,
we may assume that $\zeta\in \{0,\zeta_{en},\dots,\zeta_{en}^n\}$.
Set $i\in \{0,\dots,n\}$ such that $\zeta=\zeta_{en}^i$ or $i=0=\zeta$.
Using once again Lemma \ref{lemmadivisors} (iii) and
Lemma \ref{hyper1},  one may find an equation $X_i=\zeta'X_j$, not satisfied
by $x$, and such that $\ker(g_u-1)$ must be included in the corresponding
hyperplane; hence (iii) may not hold.
\end{proof}

The above proposition tells us about $T$-divisibility among elements
of $P_G$. Can an element of $G(e,e,n+1)$ not associated to a non-crossing
partition divide an element of $P_G$? We will prove below that this
may not happen. We start with the easier case of reflections.

\begin{lemma}
\label{goodreflections}
Let $t\in T$.
Then $t\in P_G \Leftrightarrow t \preccurlyeq_T c$.
\end{lemma}

\begin{proof}
The $\Rightarrow$ implication follows from Proposition \ref{precTdiv}.
To prove the converse, we proceed by direct computation. The reflections
not in $P_G$ have hyperplanes of equation
$\zeta X_i=X_j$, with $i,j\in \{1,\dots,n\}$,
$i< j$
and $\zeta\in \mu_e-\{1,\zeta_e\}$ (see the explicit list of reflections
in $P_G$ from Lemma \ref{hyper1}).

It is convenient to introduce condensed notation to describe
monomial matrices:
if $M$ is an $(n+1)\times (n+1)$ monomial matrix (with lines and
columns indexed by $\{0,\dots,n\}$), we represent $M$ by
$$\begin{pmatrix}
i_1 & i_2 & \dots & i_p \\
\alpha_1 & \alpha_2 & \dots & \alpha_p 
\end{pmatrix}
\begin{pmatrix}
j_1 & j_2 & \dots & j_{q} \\
\beta_1 & \beta_2 & \dots & \beta_q
\end{pmatrix}
\begin{pmatrix}
k_1 & k_2 & \dots & k_{r} \\
\gamma_1 & \gamma_2 & \dots & \gamma_r
\end{pmatrix}
\dots$$
where $$\begin{pmatrix}
i_1 & i_2 & \dots & i_p 
\end{pmatrix}
\begin{pmatrix}
j_1 & j_2 & \dots & j_{q} 
\end{pmatrix}
\begin{pmatrix}
k_1 & k_2 & \dots & k_{r} 
\end{pmatrix}
\dots$$
is the cycle decomposition of the underlying permutation $\sigma_M$
(including singleton cycles),
and the second line gives corresponding non-trivial coefficient:
for example, if $n=3$, by
$$
\begin{pmatrix}
0 \\ \zeta \end{pmatrix}
\begin{pmatrix}
1 & 2 & 3\\
\lambda & \mu & \nu
\end{pmatrix},$$
we mean the matrix
$$\begin{pmatrix}
\zeta & 0 & 0 & 0\\
0 & 0 & 0 & \nu \\
0 & \lambda & 0 & 0\\
0 & 0 & \mu & 0 
\end{pmatrix}$$

The order of the factors does not matter. Each factor is called
a \emph{generalised cycle}.

With this notation system,
$c$ is represented
by
$$\begin{pmatrix}
0  \\
\zeta_e^{-1}
\end{pmatrix}
\begin{pmatrix}
n & n-1 & \dots & 2 & 1 \\
1 & 1 & \dots & 1 & \zeta_e
\end{pmatrix}$$
and the reflection $t$ in $G(e,e,n+1)$ with hyperplane $X_i=\zeta X_j$
($1\leq i < j \leq n$) is
represented by
$$\begin{pmatrix}
i & j \\
\zeta & \zeta^{-1}
\end{pmatrix}$$
The rule for multiplying such symbols is easy to figure out. The product
$tc$ is represented by
$$\begin{pmatrix}
0  \\
\zeta_e^{-1}
\end{pmatrix}
\begin{pmatrix}
n & n-1 & \dots & j+2 & j+1        & i  & i-1 & \dots & 1\\
1 & 1 & \dots   & 1   & \zeta^{-1} & 1  & 1   & \dots & \zeta_e
\end{pmatrix}
\begin{pmatrix}
j & j-1 & \dots & i+2 & i+1 \\
1 & 1 & \dots   &   1 & \zeta
\end{pmatrix}
$$
We have to prove that if $\zeta\notin \{1,\zeta_e\}$, the kernel of $tc-1$
is trivial. This follows from the elementary remark: a monomial
matrix $M\in\GL(\BC^{n+1})$ has some non-trivial fixed points if and only
if it has a generalised factor whose product of coefficients is $1$
(actually, the number of such factors is $\dim (\ker(M-1))$) --
in our situation, we have three generalised cycles with product of
coefficients $\zeta_e^{-1}$, $\zeta^{-1}\zeta_e$ and $\zeta$.
\end{proof}

\begin{proposition}
\label{divisors}
For all $g \in G(e,e,n+1)$, we have
$$g\in P_G \; \; \Leftrightarrow \; \; 
g \preccurlyeq_T c.$$
\end{proposition}

In the proof (and later in the text), we use the notion of \emph{parabolic
subgroup} of a reflection group: if $W\subseteq \GL(V)$,
if $X\subseteq V$, the associated parabolic subgroup
$W_X$ is $\{w\in W | \forall x\in X, w(x)=x\}$.

\begin{proof}
The implication
$\Rightarrow$ follows from Proposition \ref{precTdiv}, (i) $\Rightarrow$ (ii),
and from the fact that $c=g_{\triv}$.

We prove the converse implication by induction on $n$.
Let $g\in G(e,e,n+1)$ such that $g\preccurlyeq_T c$. Let $t \in T$ such
that $t\preccurlyeq_T g$. Set $g':=t^{-1}g$, $c':=t^{-1}c$.
Since $t \preccurlyeq_T g \preccurlyeq_T c$, we deduce from Lemma
\ref{goodreflections} that $t\in P_G$. In particular, there is a height $1$
non-crossing partition $u$ such that $t=g_u$, and 
$c'=g_{\overline{u}}$.
By Proposition \ref{precTdiv}, only reflections
in the parabolic
subgroup $G'$ fixing $\ker(c'-1)$ may appear in reduced $T$-decompositions
of $c'$; let $T'$ be the set of reflections in $G'$.
We discuss by cases:
\begin{itemize}
\item If $t'$ is (short) symmetric, $\overline{u}$ is long symmetric.
Denote by $\overline{u}_1$ the long symmetric
element of $NCP(e,e,n+1)$ with
only non-trivial part this long part. Denote by $\overline{u}_2$
the (short symmetric) element of $NCP(e,e,n+1)$ with non-singleton parts
the remaining non-singleton parts of $\overline{u}$.
\begin{figure}[ht]
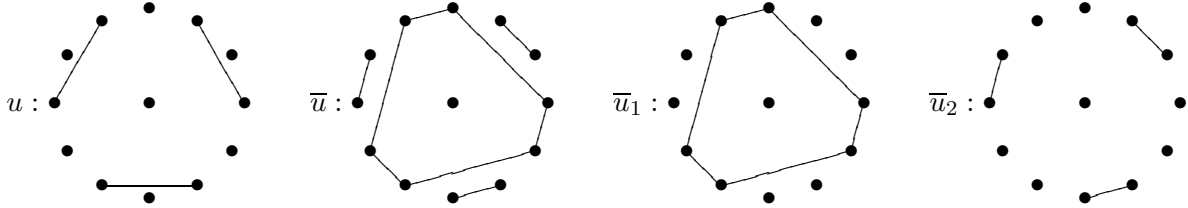

$$
u: \xy/r3pc/:,{\xypolygon12"A"{~={0}~>{}}},
"A1";"A3"**@{-},"A5";"A7"**@{-},"A9";"A11"**@{-},
"A0"*{\bullet},"A1"*{\bullet},"A2"*{\bullet},"A3"*{\bullet},"A4"*{\bullet},
"A5"*{\bullet},"A6"*{\bullet},"A7"*{\bullet},"A8"*{\bullet},
"A9"*{\bullet},"A10"*{\bullet},"A11"*{\bullet},"A12"*{\bullet}\endxy
\qquad
\overline{u}:
\xy/r3pc/:,{\xypolygon12"A"{~={0}~>{}}},
"A2";"A3"**@{-},"A6";"A7"**@{-},"A10";"A11"**@{-},
"A1";"A4"**@{-},"A4";"A5"**@{-},
"A5";"A8"**@{-},"A8";"A9"**@{-},
"A9";"A12"**@{-},"A12";"A1"**@{-},
"A0"*{\bullet},"A1"*{\bullet},"A2"*{\bullet},"A3"*{\bullet},"A4"*{\bullet},
"A5"*{\bullet},"A6"*{\bullet},"A7"*{\bullet},"A8"*{\bullet},
"A9"*{\bullet},"A10"*{\bullet},"A11"*{\bullet},"A12"*{\bullet}\endxy
\qquad
\overline{u}_1:
\xy/r3pc/:,{\xypolygon12"A"{~={0}~>{}}},
"A1";"A4"**@{-},"A4";"A5"**@{-},
"A5";"A8"**@{-},"A8";"A9"**@{-},
"A9";"A12"**@{-},"A12";"A1"**@{-},
"A0"*{\bullet},"A1"*{\bullet},"A2"*{\bullet},"A3"*{\bullet},"A4"*{\bullet},
"A5"*{\bullet},"A6"*{\bullet},"A7"*{\bullet},"A8"*{\bullet},
"A9"*{\bullet},"A10"*{\bullet},"A11"*{\bullet},"A12"*{\bullet}\endxy
\qquad
\overline{u}_2:
\xy/r3pc/:,{\xypolygon12"A"{~={0}~>{}}},
"A2";"A3"**@{-},"A6";"A7"**@{-},"A10";"A11"**@{-},
"A0"*{\bullet},"A1"*{\bullet},"A2"*{\bullet},"A3"*{\bullet},"A4"*{\bullet},
"A5"*{\bullet},"A6"*{\bullet},"A7"*{\bullet},"A8"*{\bullet},
"A9"*{\bullet},"A10"*{\bullet},"A11"*{\bullet},"A12"*{\bullet}\endxy
$$
\caption{Decomposing the partition $\overline{u}$, when $u$ is symmetric}
\end{figure}
The elements $c'_1:=g_{\overline{u}_1}$ and $c'_2:=g_{\overline{u}_2}$
satisfy $c'_1c'_2 = c'_2 c'_1 = c'$. They are Coxeter elements in 
parabolic subgroups $G'_1$ and $G'_2$
of respective types $G(e,e,p)$ and $G(1,1,q)$, where
$p$ and $q$ satisfy $p+q=n$.
Set $A'_1:=\{t\in T | t\preccurlyeq_T c'_1\}$,
$A'_2:=\{t\in T | t\preccurlyeq_T c'_2\}$ and
$A':=\{t\in T | t\preccurlyeq_T c'\}$.
Elements of $A'_1$ (resp. $A'_2$)
correspond to height $1$ non-crossing
partitions below $\overline{u}_1$ (resp. $\overline{u}_2$),
and $A' = \pi(A) \cap W' = A'_1 \cup A'_2$.
Since $g'\preccurlyeq c'$, 
a reduced $T$-decomposition of $g'$ consists of elements of $A'$.
The elements of $A'_1$ commute with the elements of $A'_2$; by regrouping
them, we obtain a decomposition $g'=g'_1g'_2$, with $g'_1\in G'_1$ and
$g'_2\in G'_2$ satifying $g'_1\preccurlyeq_{T'_1} c'_1$ and
$g'_2\preccurlyeq_{T'_2} c'_2$, where $T'_1 = T \cap G'_1$ and
$T'_2 = T\cap G'_2$.
By induction assumption, $g'_1=g_v$, where $v$ belongs to $NCP(e,e,p)$
(or, more rigorously, to the image of $NCP(e,e,p)$ under the identification
of $NCP(e,e,p)$ with elements of $NCP(e,e,n+1)$ below $\overline{u}_1$).
Using the similar result (from \cite{dualmonoid}) for the type $A$ case,
we also know that $g'_2 = g_w$, where $w$ belongs to $NCP(1,1,q)$
(or, more rigorously, to the image of $NCP(1,1,q)$ under the identification
of $NCP(1,1,q)$ with elements of $NCP(e,e,n+1)$ below $\overline{u}_2$;
the identification sends each non-singleton part to $e$ copies, one for each
element of the orbit for the action of $\mu_e$ on the set of parts of
$\overline{u}_2$). Set $u' := v\vee w$. We have $g'=g_{u'}$.
The partition $u\vee u'$ is easy to construct, and one may observe
that $u \backslash (u\vee u')= u'$. This implies that
$g_{u\vee u'} = g_u g_{u'} = t g' = g$.
We have proved our claim that $g\in P_G$.
\item When $u$ is asymmetric, the parabolic subgroup $G'$ is of type
$G(1,1,n+1)$.
\begin{figure}[ht]
$$
u: \xy/r3pc/:,{\xypolygon12"A"{~={0}~>{}}},
"A0";"A3"**@{-},
"A0"*{\bullet},"A1"*{\bullet},"A2"*{\bullet},"A3"*{\bullet},"A4"*{\bullet},
"A5"*{\bullet},"A6"*{\bullet},"A7"*{\bullet},"A8"*{\bullet},
"A9"*{\bullet},"A10"*{\bullet},"A11"*{\bullet},"A12"*{\bullet}\endxy
\qquad
\overline{u}:
\xy/r3pc/:,{\xypolygon12"A"{~={0}~>{}}},
"A0";"A4"**@{-},"A4";"A5"**@{-},"A5";"A6"**@{-},
"A6";"A7"**@{-},"A7";"A0"**@{-},
"A8";"A9"**@{-},"A9";"A10"**@{-},
"A10";"A11"**@{-},"A11";"A8"**@{-},
"A12";"A1"**@{-},"A1";"A2"**@{-},
"A2";"A3"**@{-},"A3";"A12"**@{-},
"A0"*{\bullet},"A1"*{\bullet},"A2"*{\bullet},"A3"*{\bullet},"A4"*{\bullet},
"A5"*{\bullet},"A6"*{\bullet},"A7"*{\bullet},"A8"*{\bullet},
"A9"*{\bullet},"A10"*{\bullet},"A11"*{\bullet},"A12"*{\bullet}\endxy$$
\caption{The element $\overline{u}$ when $u$ is asymmetric.}
\end{figure}
We conclude with a discussion similar (but simpler) to the above
one (the identification of $NCP(1,1,n+1)$ with the interval
$[1,\overline{u}]$ in 
$NCP(e,e,n+1)$ is Lemma \ref{intervals} (ii)).
\end{itemize}
\end{proof}

\begin{lemma}
Let $g\in P_G$. Choose $(t_1,\dots,t_k)\in \Red_T(g)$. Consider
the non-crossing partitions $u_1,\dots,u_k$ of height $1$ 
corresponding to $t_1,\dots,t_k$.
Then $g_{u_1\vee \dots \vee u_k} = g$.
The element $u_1\vee \dots \vee u_k$ only depends
on $g$, and not on the choice
of $(t_1,\dots,t_k)$ in $Red_T(g)$.
\end{lemma}

\begin{proof}
By assumption, there is an element $w\in NCP(e,e,n+1)$ such
that $g_w=g$.
By \ref{lemmadivisors} (ii), we have
$$\codim (\ker (g-1)) = l_T(g) = k.$$
Since $g=t_1\dots t_k$, the only possibility is that the reflecting 
hyperplanes of $t_1\dots t_k$ intersect transversally and
$\ker(g-1) = \bigcap_{i=1}^k \ker(t_i-1)$.

Set $u:=u_1\vee \dots \vee u_k$.
By \ref{lemmadivisors} (iii),
we have $\ker(g_u-1) \subseteq \bigcap_{i=1}^k \ker(t_i-1)=\ker(g_w-1)$ thus,
by Proposition \ref{precTdiv}, $w \preccurlyeq u$.
On the other hand, for
all $i$, $\ker(g_w-1) \subseteq \ker(t_i-1)=\ker(g_{u_i}-1)$, thus
$u_i\preccurlyeq w$.
We conclude, using the definition of the ``meet'' operation, that $u=w$.
\end{proof}

With the notations of the lemma, set $u_g:=u_1\vee \dots \vee u_k$.
The lemma precisely
asserts that $g\mapsto u_g$ is an inverse of the
map $NCP(e,e,n+1)\rightarrow G(e,e,n+1), u\mapsto g_u$.
It is then clear that we have the following theorem, which summarises
the results of this section:

\begin{theorem}
\label{titssection}
The diagram
$$\xymatrix{
 & P_B \ar[dd]^{\pi}
 \\
NCP(e,e,n+1) \ar[ur]^{u\mapsto b_u} \\
& P_G \ar[ul]^{g\mapsto u_g} }$$
is commutative. Its arrows are poset isomorphisms, where
\begin{itemize}
\item $NCP(e,e,n+1)$
is endowed with the ``is finer than'' relation $\preccurlyeq$,\\
\item $P_B$ is endowed with the relation 
$\preccurlyeq_{P_B}$
defined by $\forall b,b'\in P_B,
b' \preccurlyeq_{P_B} b \Leftrightarrow \exists b''\in P_B,
b'b''=b$, \\
\item $P_G$ is endowed with the ``is initial segment
of reduced $T$-decomposition'' relation $\preccurlyeq_T$.
\end{itemize}
\end{theorem}

In the above theorem, the composed map $g\mapsto b_{u_g}$ is an analog of
the ``Tits section'' from a finite Coxeter group to its Artin group.

\section{A Garside structure for $B(e,e,n+1)$}
\label{section5}

\begin{definition}
The \emph{dual braid monoid} of type $(e,e,n+1)$ is the
submonoid $M(e,e,n+1)$ of $B(e,e,n+1)$ generated by $P_B$.
\end{definition}

We prove in this section that $M(e,e,n+1)$ is Garside monoid, with Garside
element $b_{\triv}$ and set of simples $P_B$. The proof uses the results
from the previous section, about reduced decompositions in $G(e,e,n+1)$,
and is very similar to the proof in \cite{bdm}.

Let us complete the commutative diagram of Theorem \ref{titssection}
by adding the natural inclusions:

$$\xymatrix{
 & P_B \ar[dd]^{\pi}
 \ar@{^{(}->}[r] & M(e,e,n+1) \ar@{^{(}->}[r] & B(e,e,n+1) \ar@{>>}[dd]^{\pi}
 \\
NCP(e,e,n+1) \ar[ur] \\
& P_G \ar[ul] \ar@{^{(}->}[rr] & & G(e,e,n+1) }$$

\begin{theorem}
\label{maintheorem}
The monoid $M(e,e,n+1)$ is a Garside monoid. The element
$b_{\triv}$ is a Garside element, for which the set of simple elements
is $P_B$. The atoms of $M(e,e,n+1)$ are the images by $\alpha$ of 
height $1$ non-crossing partitions.

The inclusion $M(e,e,n+1) \hookrightarrow B(e,e,n+1)$ identifies 
$B(e,e,n+1)$ with the group of fractions of $M(e,e,n+1)$. Hence
$B(e,e,n+1)$ is a Garside group.
\end{theorem}

For general properties of Garside groups, a good reference is \cite{dehgar}.

\begin{proof}
Since $P_G = \{g \in G(e,e,n+1) | g \preccurlyeq_T c\}$
(Proposition \ref{divisors}), it may
be endowed with the partial product structure defined in section
0.4 of \cite{dualmonoid} (note that, since $T$ is
stable under conjugacy, the relations
$\prec_T$ and $\succ_T$ of \emph{loc. cit.} coincide -- in particular,
$c$ is balanced).
By Theorem \ref{titssection}, $(P_G,\preccurlyeq_T)$ is isomorphic
to $(NCP(e,e,n+1),\preccurlyeq)$, and in particular
(Lemma \ref{latticeproperty}) it is a lattice.
This implies, using \cite[Theorem 0.5.2]{dualmonoid}, that the monoid
$\mathbf{M}(P_G)$ generated by $P_G$ is a Garside monoid
(with Garside element corresponding to $c$, and set of simples
corresponding to $P_G$).

The bijection $P_G\simeq P_B$ of Theorem \ref{titssection}
is compatible with the partial product, it induces a monoid
morphism $\phi:\mathbf{M}(P_G) \rightarrow M(e,e,n+1)$, which extends
to a group morphism $\psi:\mathbf{G}(P_G) \rightarrow G(e,e,n+1)$.
We have $\phi(c) = b_{\triv}$.
The point is to prove that $\phi$ and $\psi$ are isomorphisms.
By definition, $P_B$ generates
$M(e,e,n+1)$, thus $\phi$ is surjective. By Theorem \ref{bmr2},
$P_B$ contains group generators for $B(e,e,n+1)$, so $\psi$ is also 
surjective. Since $\mathbf{M}(P_G)$ is a Garside
monoid, the canonical map $\mathbf{M}(P_G) \rightarrow \mathbf{G}(P_G)$
is injective. Thus, to complete the proof, it is sufficient to prove
that $\psi$ is injective.

To prove the injectivity of $\psi$, it is sufficient to check 
that the Brou\'e-Malle-Rouquier relations between the specific
generators in Theorem \ref{bmr2} are consequences of the
relations in $\mathbf{G}(P_G)$. We illustrate this straightforward computation
with the relation
$\langle \tau_2 \tau'_2 \rangle^e = \langle \tau'_2 \tau_2 \rangle^e$.
The generator $\tau_2$ (resp. $\tau_2'$) is identified with $a_0$ (resp.
$a_n$). The least upper bound of the partitions $u_0$ and
$u_n$ is the height $2$ partition $v:=((u_0)^{\sharp})_*=
((u_n)^{\sharp})_*$, whose only non-singleton part is $\mu_e\cup\{0\}$.
By definition of the complement operation, we have $u_0 \backslash
v= u_n$. By Proposition \ref{mappingb}, we then have
$a_0a_n = b_v$. This relation holds in $P_B$, thus is $P_G$
(via the isomorphism of Theorem \ref{titssection}) and in $\mathbf{G}(P_G)$.
By symmetry, we have relations $a_{kn}a_{(k+1)n} = b_v$, for all
integer $k$.
If $e$ is even, we have  
$$\langle a_0^{(0)}a_0^{(1)}\rangle^e =
(a_0a_n)^{e/2} = b_v^{e/2}= a_na_{2n} a_{3n}a_{4n} 
\dots a_{(e-1)n}a_{en} = a_n b_v^{e/2-1} a_0 =
\langle a_na_0\rangle^e$$
as a consequence of relations of $\mathbf{G}(P_G)$. The $e$ odd case,
as well as the other relations of Theorem \ref{bmr1},
may be obtained
in a similar manner.
\end{proof}

{\flushleft \bf An application to centralisers of periodic elements.}
The main result in \cite{bdm} is a description of centralisers of
roots of central elements in type $A$ braid groups. Similarly, we obtain
the following result:

\begin{proposition}
The element $b_{\triv}^{\frac{en}{e\wedge (n+1)}}$ generates the center 
of $B(e,e,n+1)$.
Let $k\in \BZ{\geq 1}$,
and let $k':=\frac{en}{k(n+1) \wedge en}$.
Assume $k'\neq 1$.
Then the centraliser $C_{B(e,e,n+1)}(b_{\triv}^k)$ is isomorphic
to $B(k',1,en/k')$.
\end{proposition}

\begin{proof}
The first statement, about $b_{\triv}^{\frac{en}{e\wedge (n+1)}}$,
may already be found in \cite{bmr}.

Let $k\in \BZ$.

Let $b\in P_B$. There is a unique $u\in NCP(e,e,n+1)$ such that
$b=b_u$. Since $b_ub_{\overline{u}}= b_{\triv} = b_{\overline{u}}
b_{\overline{\overline{u}}}$, we have
$$b_{\overline{\overline{u}}} = b_{\triv}^{-1} b_u b_{\triv}.
$$
Replacing $\overline{\overline{u}}$ by its description from
Lemma \ref{lemmacon} (v), we obtain
$
b_{\zeta_{en}^{n+1} u} = b_{\triv}^{-1} b_u b_{\triv}
$
and
\begin{equation}
\label{equZ}
b_{\zeta_{en}^{k(n+1)} u} = b_{\triv}^{-k} b_u b_{\triv}^k.
\end{equation}

Let $b\in C_{B(e,e,n+1)}(b_{\triv}^k)$. Write $b$ in Garside normal
form: $b = b_{\triv}^N b_1\dots b_m$, with $N\in \BZ$ and
$b_i\in P_B$. By (\ref{equZ}), conjugating by a power of $b_{\triv}$ preserves
$P_B$. If $b$ is invariant by such a conjugacy, the unicity of the normal
form implies that each $b_i$ is invariant by the conjugacy, that is,
using (\ref{equZ}), we have 
$\zeta_{en}^{k(n+1)} u_i=u_i$, where $u_i$ is the non-crossing partition
such that $b_i=b_{u_i}$.

The integer $k'$ is the order of
$\zeta':=\zeta_{en}^{k(n+1)}$. If $k'\neq 1$,
multiplication by $\zeta'$ fixes no asymmetric partition.
An symmetric partition is fixed if and only if it lies
in $NCP(k',1,en/k')$.
This implies that the natural morphism $B(k',1,en/k') \rightarrow
C_{B(e,e,n+1)}(b_{\triv}^k)$ is surjective.

To prove the injectivity, one proceeds as follows.
By \cite{bdm}, Proposition 3.26,
the monoid $C_{M(e,e,n+1)}(b_{\triv}^k)$ is a Garside monoid,
with set of simples in bijection with $NCP(k',1,en/k')$.
There is a ``dual monoid'' $M(k',1,en/k')$ for $B(k',1,en/k')$,
with $NCP(k',1,en/k')$. The bijection between their set of simples
induces an isomorphism $M(k',1,en/k') \simeq C_{M(e,e,n+1)}(b_{\triv}^k)$.
\end{proof}

\section{Transitivity of Hurwitz action}
\label{section6}

This section contains complements about the structure of $\Red_T(c)$,
which will be needed in the next section to write a simple presentation
for $G(e,e,n+1)$.

Let $g\in G(e,e,n+1)$, with $l_T(g)=k$. 
For any $i=1,\dots,k-1$ and any $(t_1,\dots,t_k)\in\Red_T(g)$,
set
$$\sigma_i ((t_1,\dots,t_k)):=(t_1,\dots,t_{i-1},t_{i}t_{i+1}t_i^{-1},t_i,
t_{i+2},\dots,t_k).$$
The right-hand side clearly belongs to $\Red_T(g)$, and one checks
that $\sigma_1,\dots,
\sigma_{k-1}$ satisfy the defining the classical braid relations
defining the Artin group of type $A_{k-1}$. The corresponding action
of this Artin group on $\Red_T(g)$ is called \emph{Hurwitz action}.
The following result extends \cite{dualmonoid}, Proposition 1.6.1.

\begin{proposition}
\label{matsumoto}
Let $g\in P_G$. 
The Hurwitz action is transitive
on $\Red_T(g)$.
\end{proposition}

We need two lemmas:

\begin{lemma}
\label{orbits}
Let $w$ be a maximal short symmetric element of $NCP(e,e,n+1)$.
We have $ht(w)=n-1=l_T(g_w)$.
There exists an element $(t_1,\dots,t_{n-1})\in\Red_T(g_w)$ such
that for all symmetric $t\in T \cap P_G$,
one may find
$i\in \{1,\dots,n-1\}$ and $\zeta\in\mu_{en}$ such that $t=\zeta t_i$.
\end{lemma}

\begin{proof}
The statement about $ht(w)$ is easy. Multiplication by $\mu_{en}$ is
transitive on the set of maximal short symmetric elements,
so it suffices to prove the lemma for a particular one.

Consider the short symmetric
non-crossing partition $u$ whose non-singleton parts are
$$\{1,\zeta_{en}\},\{\zeta_{en}^{-1},\zeta_{en}^2\},
\dots,\{\zeta_{en}^{-i+1},\zeta_{en}^i\},\dots,
\{\zeta_{en}^{-\lceil \frac{n-1}{2}\rceil+1},
\zeta_{en}^{\lceil \frac{n-1}{2}\rceil}\}$$
and their images under multiplication by $\mu_e$;
consider the short symmetric
non-crossing partition $v$ whose non-singleton parts are
$$\{1,\zeta_{en}^2\},\{\zeta_{en}^{-1},\zeta_{en}^{3}\},
\dots,\{\zeta_{en}^{-i+1},\zeta_{en}^{i+1}\},\dots,
\{\zeta_{en}^{-\lfloor \frac{n-1}{2}\rfloor+1},
\zeta_{en}^{\lfloor \frac{n-1}{2}\rfloor+1}\}$$
and their images under multiplication by $\mu_e$.
For all integer $p\geq 0$, set $s_p:=t_{-p,p+1}$ and
$s'_p:=t_{2p+1,-2p-1}$.
The element $u\vee v$ is a maximal short symmetric non-crossing partition,
with parts
$$\{0\},\{\zeta_{en}^{-\lceil \frac{n-1}{2}\rceil+1},
\zeta_{en}^{-\lceil \frac{n-1}{2}\rceil+2},\dots,
\zeta_{en}^{\lfloor \frac{n-1}{2}\rfloor+1}\}$$
and their images under multiplication by $\mu_e$.
\begin{figure}[ht]
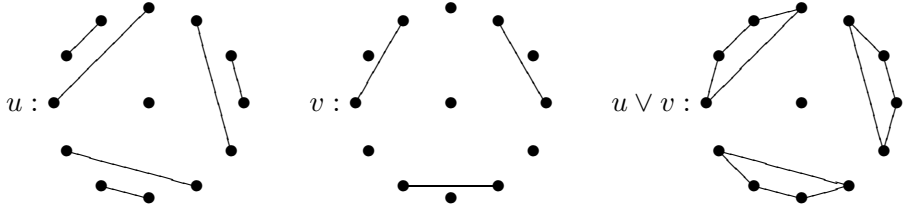

$$
u: \xy/r3pc/:,{\xypolygon12"A"{~={0}~>{}}},
"A1";"A2"**@{-},"A5";"A6"**@{-},"A9";"A10"**@{-},
"A12";"A3"**@{-},"A4";"A7"**@{-},"A8";"A11"**@{-},
"A0"*{\bullet},"A1"*{\bullet},"A2"*{\bullet},"A3"*{\bullet},"A4"*{\bullet},
"A5"*{\bullet},"A6"*{\bullet},"A7"*{\bullet},"A8"*{\bullet},
"A9"*{\bullet},"A10"*{\bullet},"A11"*{\bullet},"A12"*{\bullet}\endxy
\qquad
v:
\xy/r3pc/:,{\xypolygon12"A"{~={0}~>{}}},
"A1";"A3"**@{-},"A5";"A7"**@{-},"A9";"A11"**@{-},
"A0"*{\bullet},"A1"*{\bullet},"A2"*{\bullet},"A3"*{\bullet},"A4"*{\bullet},
"A5"*{\bullet},"A6"*{\bullet},"A7"*{\bullet},"A8"*{\bullet},
"A9"*{\bullet},"A10"*{\bullet},"A11"*{\bullet},"A12"*{\bullet}\endxy
\qquad
u\vee v:
\xy/r3pc/:,{\xypolygon12"A"{~={0}~>{}}},
"A1";"A2"**@{-},"A2";"A3"**@{-},"A3";"A12"**@{-},"A1";"A12"**@{-},
"A5";"A6"**@{-},"A6";"A7"**@{-},"A7";"A4"**@{-},"A5";"A4"**@{-},
"A9";"A10"**@{-},"A10";"A11"**@{-},"A11";"A8"**@{-},"A9";"A8"**@{-},
"A0"*{\bullet},"A1"*{\bullet},"A2"*{\bullet},"A3"*{\bullet},"A4"*{\bullet},
"A5"*{\bullet},"A6"*{\bullet},"A7"*{\bullet},"A8"*{\bullet},
"A9"*{\bullet},"A10"*{\bullet},"A11"*{\bullet},"A12"*{\bullet}\endxy
$$
\caption{Illustration in  $NCP(3,3,5)$}
\end{figure}

Set $U:=\{w\in NCP(e,e,n+1) | w\preccurlyeq u, ht(w)=1\}$ and
 $V:=\{w\in NCP(e,e,n+1) | w\preccurlyeq v, ht(w)=1\}$
The lemma follows from the following simple observations:
\begin{itemize}
\item $ht(u)=\lceil (n-1)/2 \rceil = |U|$ and
$ht(v)=\lfloor (n-1)/2 \rfloor = |V|$;
\item any reduced decomposition of $g_u$ 
(resp. $g_v$) involves all the reflections associated to elements
of $U$ (resp. $V$);
\item $v=u\backslash (u\vee v)$, thus $g_ug_v=g_{u\vee v}$;
\item $u\vee v$ is a maximal short symmetric element which admits
a reduced $T$-decomposition involving all elements in $U\cup V$; 
\item any height $1$ symmetric non-crossing partitions is obtained
by multiplying some element of $U\cup V$ by some element in $\mu_{en}$.
\end{itemize}
\end{proof}

\begin{lemma}
\label{orbits2}
Let $w'$ be a minimal long symmetric element of $NCP(e,e,n+1)$.
We have $ht(w')=2=l_T(g_w')$.
Let $u_p$ be a height $1$ asymmetric element finer than $w'$.
For all integers $k$, $u_{p+kn}\preccurlyeq w'$, and
$$\Red_{T}(g_{w'})=\{(g_{u_p},g_{u_{p+n}}),
(g_{u_{p+n}},g_{u_{p+2n}}),\dots,(g_{u_{p+(e-1)n}},g_{u_{p}})\}.$$
\end{lemma}

\begin{proof}
The partition $w'$ has only one non-singleton part,
of the form $\{0,\zeta,\zeta_e\zeta,\zeta_e^2\zeta,\dots,\zeta_e^{e-1}\zeta\}$.
One observes that height $1$ elements finer than $w'$ are the asymmetric
partitions $v_i$ with non-singleton part $\{0,\zeta_e^i\zeta\}$.
By Proposition \ref{precTdiv} (i) $\Leftrightarrow$ (ii), all sequences
in $\Red_T(w')$ have their terms among the $g_{v_i}$'s.
One has $v_i\backslash w' = v_{i+1}$. The result follows.
\end{proof}

\begin{proof}[Proof of the proposition]
We prove the proposition by induction. Assume it to be known for all
values of $n$ smaller than the considered one.

Consider first the particular case when $g=c$.
Choose $w$ maximal short symmetric in $NCP(e,e,n+1)$.
Set $w':=w\backslash \triv$. The element $w'$ is minimal long symmetric.
Choose $(t_0,\dots,t_{n-2})\in \Red_T(w)$ satisfying the condition of
Lemma \ref{orbits}.
Let $(t_{n-1},t_{n})\in \Red_T(w')$.
We have $(t_0,\dots,t_{n})\in \Red_T(c)$.
Let $(t'_0,\dots,t'_n)$ be another element of $\Red_T(c)$.

\begin{quote}
Claim. \emph{We may find $(t_0'',\dots,t_n'')$ in the Hurwitz orbit
of $(t_0,\dots,t_n)$ such that $t_0''=t_0'$.}
\end{quote}

Proof of the claim. One first remarks that it suffices to find
$(t_0'',\dots,t_n'')$ with $t_i''=t_0'$, since one may use Hurwitz to
``slide'' $t_i''$ to the beginning of the sequence.
A standard calculation shows that
the action of  $\beta:=(\sigma_1\dots \sigma_{n})^{n+1}$ on
$\Red_T(c)$ coincides with term-by-term conjugacy by $c$, which, on the
corresponding partitions, coincides with $u\mapsto \overline{\overline{u}}$,
which itself coincides with the multiplication by $\zeta_{en}^{n+1}$
(Lemma \ref{lemmacon} (v)).
By Proposition \ref{precTdiv}, $t_0'\in P_G$.
If $t_0'$ is associated to a symmetric height $1$
partition, then, using Lemma \ref{orbits}, we may apply a suitable
power of $\beta$ to find in
the Hurwitz orbit of $(t_0,\dots,t_n)$
a sequence involving $t_0'$ (on symmetric partitions, multiplying by
$\zeta_{en}^{n+1}$ is the same as multiplying by $\zeta_{en}$).
It $t_0'$ is asymmetric, the underlying partition is obtained from
the underlying partition of $t_{n-1}$ by multiplying by some $\zeta\in\mu_{en}$.
Use $\sigma_n$ to multiply by $\zeta_e$ the underlying partitions of $t_{n-1}$
and $t_n$
(Lemma \ref{orbits2}); use $\beta$ to multiply by $\zeta_{en}^{n+1}$;
to conclude, observe that $\zeta_e$ and
$\zeta_{en}^{n+1}=\zeta_e\zeta_{en}$ generate $\mu_{en}$.

To conclude that $(t_0,\dots,t_n)$ and $(t_0',\dots,t_n')$ lie
in the same Hurwitz orbit, it is enough to prove that
$(t_1',\dots,t_n')$ and $(t_1'',\dots,t_n'')$ lie in the same
Hurwitz orbit. This follows either from the induction assumption
(when the parabolic subgroup fixing $\ker(t_0'c-1)$ is of type
$G(e,e,n)$) or from \cite{dualmonoid} 1.6.1 (when this parabolic
subgroup is a real reflection group).

Similarly, the general case follows either from the induction
assumption or from \cite{dualmonoid} 1.6.1, depending on the type
of the parabolic subgroup $W'$ fixing $\ker(g-1)$. This is because
$\Red_T(g)$ coincides with $\Red_{T'}(g)$, where $T'=T\cap W'$
(Proposition \ref{precTdiv}).
\end{proof}

\section{Explicit presentations of $M(e,e,n+1)$ and $B(e,e,n+1)$}
\label{section7}

The generating set, in the theorem below, is the set $A$ of
elements associated with non-crossing partitions of height $1$
(section \ref{section3}). By Theorem \ref{titssection}, it is
in bijection with $T\cap P_c$.
Alternatively, we have $A=\{ a\in P_B | l(a)=1\}$ (where $l$ is the
canonical length function on $B(e,e,n+1)$).
The elements of $A$ were explicitely described in section \ref{section3}:
there are $n(n-1)$ symmetric elements $a_{p,q}$, and $en$ asymmetric
elements $a_r$ (where $p,q,r$ are arbitrary integers such that $|p-q|<n$;
recall that $a_{p,q}=a_{q,p}= a_{p+n,q+n}$ and $a_{p}=a_{p+en}$).

\begin{theorem}
\label{listofrels}
(1) The following length $2$ relations between elements of $A$ hold
in $M(e,e,n+1)$:

For every quadruple $(p,q,r,s)$ of integers such that
$p<q<r< s < p+n$, there are relations:
$$\begin{array}{ lcl}
 (\R_1) &  a_{p,q} a_{r,s} = a_{r,s} a_{p,q}    &   \\
 (\R_2) &  a_{p,s} a_{q,r} = a_{q,r} a_{p,s} .   &     
\end{array}$$
For every triple $(p,q,r)$ of integers such that
$p<q<r< p+n$, there are relations:
$$\begin{array}{lcl}
  (\R_3) &   a_{p,r} a_{q,r}  =  a_{q,r} a_{p,q} =  a_{p,q} a_{p,r}  &   \\
  (\R_4) &     a_{p,q} a_r =    a_r   a_{p,q} 
.                  &
\end{array}$$
For every pair $(p,q)$ of integers such that $p<q<p+n$,
there are relations:
$$\begin{array}{lcl}
  (\R_5) &   a_{p,q} a_{p}  =  a_{p} a_{q} =  
a_{q}a_{p,q}  & \\
\end{array}$$
For every $p\in \BZ/n\BZ$, there are relations:
$$ (\R_6) \quad  \quad 
a_p a_{p+n} 
= \cdots = a_{p+in} a_{p+(i+1)n}  = \cdots =  a_{p+(e-1)n} a_p.$$  

\begin{figure}[p]
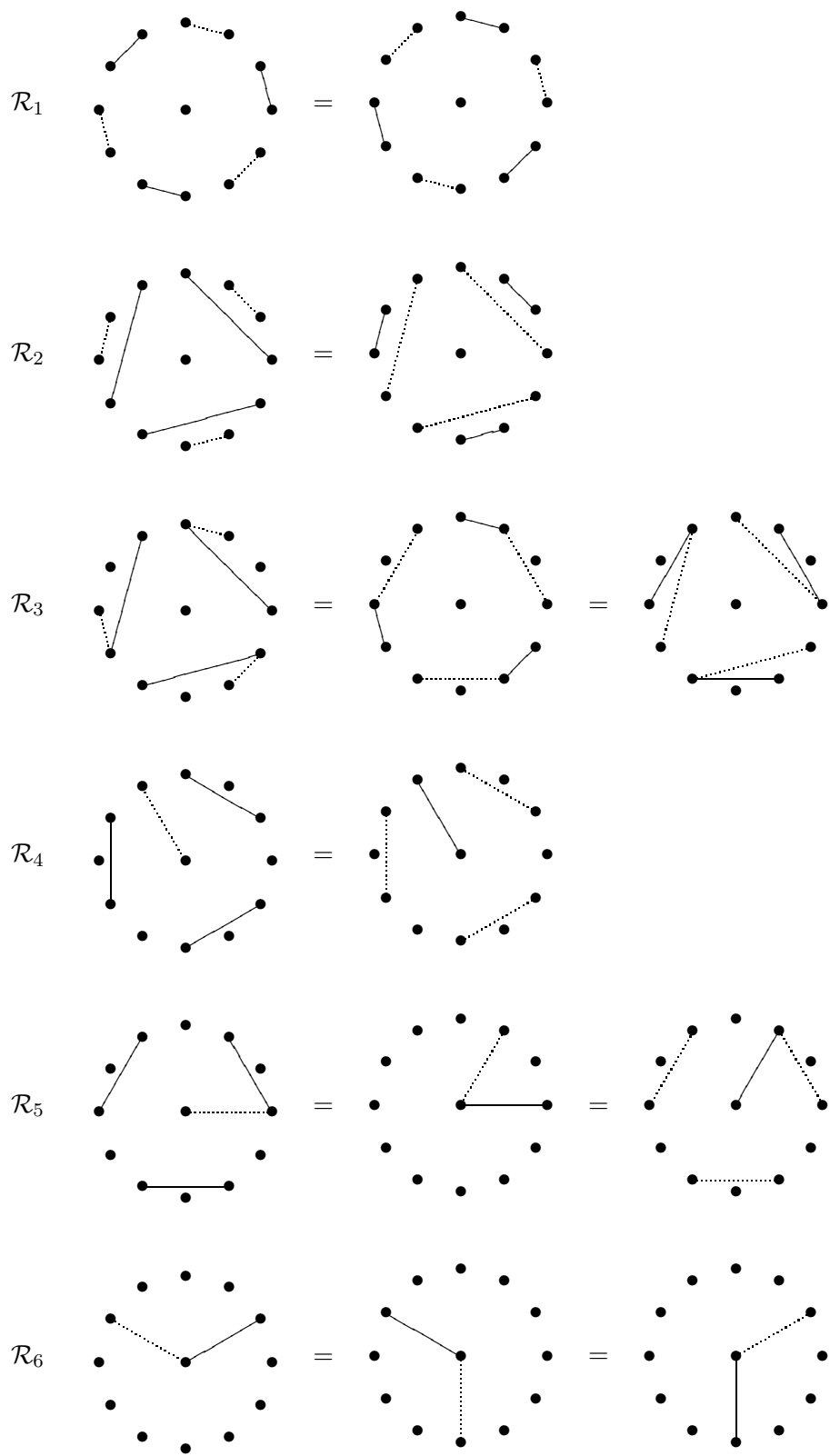

\begin{eqnarray*}
\R_1
& &
\xy/r3pc/:,{\xypolygon12"A"{~={0}~>{}}},
"A1";"A2"**@{-},"A5";"A6"**@{-},"A9";"A10"**@{-},
"A3";"A4"**@{.},"A7";"A8"**@{.},"A11";"A12"**@{.},
"A0"*{\bullet},"A1"*{\bullet},"A2"*{\bullet},"A3"*{\bullet},"A4"*{\bullet},
"A5"*{\bullet},"A6"*{\bullet},"A7"*{\bullet},"A8"*{\bullet},
"A9"*{\bullet},"A10"*{\bullet},"A11"*{\bullet},"A12"*{\bullet}\endxy
\quad = \quad
\xy/r3pc/:,{\xypolygon12"A"{~={0}~>{}}},
"A3";"A4"**@{-},"A7";"A8"**@{-},"A11";"A12"**@{-},
"A1";"A2"**@{.},"A5";"A6"**@{.},"A9";"A10"**@{.},
"A0"*{\bullet},"A1"*{\bullet},"A2"*{\bullet},"A3"*{\bullet},"A4"*{\bullet},
"A5"*{\bullet},"A6"*{\bullet},"A7"*{\bullet},"A8"*{\bullet},
"A9"*{\bullet},"A10"*{\bullet},"A11"*{\bullet},"A12"*{\bullet}\endxy
\\
\\
\\
\R_2 & &
\xy/r3pc/:,{\xypolygon12"A"{~={0}~>{}}},
"A1";"A4"**@{-},"A5";"A8"**@{-},"A9";"A12"**@{-},
"A3";"A2"**@{.},"A7";"A6"**@{.},"A11";"A10"**@{.},
"A0"*{\bullet},"A1"*{\bullet},"A2"*{\bullet},"A3"*{\bullet},"A4"*{\bullet},
"A5"*{\bullet},"A6"*{\bullet},"A7"*{\bullet},"A8"*{\bullet},
"A9"*{\bullet},"A10"*{\bullet},"A11"*{\bullet},"A12"*{\bullet}\endxy
\quad = \quad
\xy/r3pc/:,{\xypolygon12"A"{~={0}~>{}}},
"A3";"A2"**@{-},"A7";"A6"**@{-},"A11";"A10"**@{-},
"A1";"A4"**@{.},"A5";"A8"**@{.},"A9";"A12"**@{.},
"A0"*{\bullet},"A1"*{\bullet},"A2"*{\bullet},"A3"*{\bullet},"A4"*{\bullet},
"A5"*{\bullet},"A6"*{\bullet},"A7"*{\bullet},"A8"*{\bullet},
"A9"*{\bullet},"A10"*{\bullet},"A11"*{\bullet},"A12"*{\bullet}\endxy
\\
\\
\\
\R_3 & &
\xy/r3pc/:,{\xypolygon12"A"{~={0}~>{}}},
"A1";"A4"**@{-},"A5";"A8"**@{-},"A9";"A12"**@{-},
"A3";"A4"**@{.},"A7";"A8"**@{.},"A11";"A12"**@{.},
"A0"*{\bullet},"A1"*{\bullet},"A2"*{\bullet},"A3"*{\bullet},"A4"*{\bullet},
"A5"*{\bullet},"A6"*{\bullet},"A7"*{\bullet},"A8"*{\bullet},
"A9"*{\bullet},"A10"*{\bullet},"A11"*{\bullet},"A12"*{\bullet}\endxy
\quad = \quad
\xy/r3pc/:,{\xypolygon12"A"{~={0}~>{}}},
"A3";"A4"**@{-},"A7";"A8"**@{-},"A11";"A12"**@{-},
"A1";"A3"**@{.},"A5";"A7"**@{.},"A9";"A11"**@{.},
"A0"*{\bullet},"A1"*{\bullet},"A2"*{\bullet},"A3"*{\bullet},"A4"*{\bullet},
"A5"*{\bullet},"A6"*{\bullet},"A7"*{\bullet},"A8"*{\bullet},
"A9"*{\bullet},"A10"*{\bullet},"A11"*{\bullet},"A12"*{\bullet}\endxy
\quad = \quad
\xy/r3pc/:,{\xypolygon12"A"{~={0}~>{}}},
"A1";"A3"**@{-},"A5";"A7"**@{-},"A9";"A11"**@{-},
"A1";"A4"**@{.},"A5";"A8"**@{.},"A9";"A12"**@{.},
"A0"*{\bullet},"A1"*{\bullet},"A2"*{\bullet},"A3"*{\bullet},"A4"*{\bullet},
"A5"*{\bullet},"A6"*{\bullet},"A7"*{\bullet},"A8"*{\bullet},
"A9"*{\bullet},"A10"*{\bullet},"A11"*{\bullet},"A12"*{\bullet}\endxy
\\
\\
\\
\R_4
& & 
\xy/r3pc/:,{\xypolygon12"A"{~={0}~>{}}},
"A2";"A4"**@{-},"A6";"A8"**@{-},"A10";"A12"**@{-},
"A5";"A0"**@{.},
"A0"*{\bullet},"A1"*{\bullet},"A2"*{\bullet},"A3"*{\bullet},"A4"*{\bullet},
"A5"*{\bullet},"A6"*{\bullet},"A7"*{\bullet},"A8"*{\bullet},
"A9"*{\bullet},"A10"*{\bullet},"A11"*{\bullet},"A12"*{\bullet}\endxy
\quad = \quad
\xy/r3pc/:,{\xypolygon12"A"{~={0}~>{}}},
"A2";"A4"**@{.},"A6";"A8"**@{.},"A10";"A12"**@{.},
"A5";"A0"**@{-},
"A0"*{\bullet},"A1"*{\bullet},"A2"*{\bullet},"A3"*{\bullet},"A4"*{\bullet},
"A5"*{\bullet},"A6"*{\bullet},"A7"*{\bullet},"A8"*{\bullet},
"A9"*{\bullet},"A10"*{\bullet},"A11"*{\bullet},"A12"*{\bullet}\endxy
\\
\\
\\
\R_5
& &
\xy/r3pc/:,{\xypolygon12"A"{~={0}~>{}}},
"A1";"A3"**@{-},"A5";"A7"**@{-},"A9";"A11"**@{-},
"A0";"A1"**@{.},
"A0"*{\bullet},"A1"*{\bullet},"A2"*{\bullet},"A3"*{\bullet},"A4"*{\bullet},
"A5"*{\bullet},"A6"*{\bullet},"A7"*{\bullet},"A8"*{\bullet},
"A9"*{\bullet},"A10"*{\bullet},"A11"*{\bullet},"A12"*{\bullet}\endxy
\quad = \quad
\xy/r3pc/:,{\xypolygon12"A"{~={0}~>{}}},
"A0";"A1"**@{-},
"A0";"A3"**@{.},
"A0"*{\bullet},"A1"*{\bullet},"A2"*{\bullet},"A3"*{\bullet},"A4"*{\bullet},
"A5"*{\bullet},"A6"*{\bullet},"A7"*{\bullet},"A8"*{\bullet},
"A9"*{\bullet},"A10"*{\bullet},"A11"*{\bullet},"A12"*{\bullet}\endxy
\quad = \quad
\xy/r3pc/:,{\xypolygon12"A"{~={0}~>{}}},
"A0";"A3"**@{-},
"A1";"A3"**@{.},"A5";"A7"**@{.},"A9";"A11"**@{.},
"A0"*{\bullet},"A1"*{\bullet},"A2"*{\bullet},"A3"*{\bullet},"A4"*{\bullet},
"A5"*{\bullet},"A6"*{\bullet},"A7"*{\bullet},"A8"*{\bullet},
"A9"*{\bullet},"A10"*{\bullet},"A11"*{\bullet},"A12"*{\bullet}\endxy
\\
\\
\\
\R_6
& & 
\xy/r3pc/:,{\xypolygon12"A"{~={0}~>{}}},
"A0";"A2"**@{-},
"A0";"A6"**@{.},
"A0"*{\bullet},"A1"*{\bullet},"A2"*{\bullet},"A3"*{\bullet},"A4"*{\bullet},
"A5"*{\bullet},"A6"*{\bullet},"A7"*{\bullet},"A8"*{\bullet},
"A9"*{\bullet},"A10"*{\bullet},"A11"*{\bullet},"A12"*{\bullet}\endxy
\quad = \quad
\xy/r3pc/:,{\xypolygon12"A"{~={0}~>{}}},
"A0";"A6"**@{-},
"A0";"A10"**@{.},
"A0"*{\bullet},"A1"*{\bullet},"A2"*{\bullet},"A3"*{\bullet},"A4"*{\bullet},
"A5"*{\bullet},"A6"*{\bullet},"A7"*{\bullet},"A8"*{\bullet},
"A9"*{\bullet},"A10"*{\bullet},"A11"*{\bullet},"A12"*{\bullet}\endxy
\quad = \quad
\xy/r3pc/:,{\xypolygon12"A"{~={0}~>{}}},
"A0";"A10"**@{-},
"A0";"A2"**@{.},
"A0"*{\bullet},"A1"*{\bullet},"A2"*{\bullet},"A3"*{\bullet},"A4"*{\bullet},
"A5"*{\bullet},"A6"*{\bullet},"A7"*{\bullet},"A8"*{\bullet},
"A9"*{\bullet},"A10"*{\bullet},"A11"*{\bullet},"A12"*{\bullet}\endxy
\end{eqnarray*}
\caption{Diagrammatic interpretation of the relations in $G(3,3,5)$.
All relations
involves products of two generators. The full segments represent
the first term, the dotted segments represent 
the second term.}
\label{figurerels}
\end{figure}

(2) Denote by $\R$ the set of all relations considered in (1).
We have presentations:
$$M(e,e,n+1) \simeq \left< A \left| \R \right. \right>_{\text{\em Monoid}}
\qquad \text{and} \qquad
B(e,e,n+1) \simeq \left< A \left| \R \right. \right>_{\text{\em Group}}.
$$
\end{theorem}

On the graphical illustration (Figure \ref{figurerels}), one may realise that:
generators commute when the non-trivial parts (``edges'') of their
associated partitions do not interest; they
satisfy another length $2$ relation when the edges have common
endpoints; there is no relation when edges cross in their inner part.
Note also than, whenever $u,v,w$ are such that
$ht(u)=1$, $ht(v)=1$, $ht(w)=2$, $u\preccurlyeq w$ and $v=u\backslash w$,
we have a relation $b_ub_v = b_v b_{v\backslash w}$, and that all relations
are of this form. Relations are in correspondence with non-crossing partitions
of height $2$.
This allows a more intrinsic (though less explicit)
reformulation of the theorem.

\begin{proof}
The observation that the relations are of the form
$b_ub_v = b_v b_{v\backslash w}$ suffices to prove (1).

To prove (2), using standard facts about Garside structures
(see \cite{dualmonoid}, Section 0), we
deduce from Theorem \ref{maintheorem} a monoid presentation for $M(e,e,n+1)$
with:
\begin{itemize}
\item Generating set: $P_G$.
\item Relations: $g'g''=g$ whenever the relation holds in
$G(e,e,n+1)$ and $l_T(g')+l_T(g'')=l_T(g)$.
\end{itemize}
One rewrites this presentation as follows:
\begin{itemize}
\item Generating set: $P_G\cap T$.
\item Relations: $t_1\dots t_k=t_1'\dots t'_k$ whenever the relation
holds in $G$ and $(t_1,\dots ,t_k)\in\Red_T(t_1\dots t_k)$
and $(t_1',\dots, t'_k)\in\Red_T(t_1'\dots t_k')=\Red_T(t_1\dots t_k)$.
\end{itemize}
Now Proposition \ref{matsumoto} implies the above relations, for $k\geq 1$,
are consequences of the relations for $k=2$.
When considering the presentation as a group presentation, we obtain
a presentation for the group of fractions $G(e,e,n+1)$ of $M(e,e,n+1)$.
\end{proof}

The presentation has some symmetries:

\begin{proposition}
The monoid $M(e,e,n+1)$ is isomorphic to the opposed monoid
$M(e,e,n+1)^{\text{op}}$.
\end{proposition}

\begin{proof}
Consider an plane axial symmetry $\phi$ preserving the regular
$ne$-gon $\mu_{ne}$. By its action on diagrams, $\phi$ induces
an involution of the generating set $A$.
The relations of Proposition \ref{listofrels} only depend on (oriented)
incidence of diagrams. By examining their graphical interpretations,
one may observe that if $ab=cd$ is a relation in $\R$
between $a,b,c,d\in A$, then $\phi(b)\phi(a)=\phi(d)\phi(c)$ is also
in $\R$. Hence $\phi$ realizes an anti-automorphism of $M(e,e,n+1)$.
\end{proof}

Note that the above proof does not provide us with a natural antiautomorphism
of $M(e,e,n+1)$. By combining all axial symmetries and rotations of
$\mu_{ne}$, we obtain a dihedral group of order $2ne$, whose even (resp. odd)
elements are automorphisms (resp. antiautomorphisms) 
of $M(e,e,n+1)$. The situation is similar to the one in \cite{dualmonoid},
Section 2.4.

{\bf \flushleft Remark.} Our presentation
also makes it easy to identify $B(e,1,n)$ with
a subgroup of $B(e,e,n+1)$.
For all $p\in \BZ/n\BZ$, set $b_p:=a_{p}a_{p+n}$  (by $\R_6$,
we have $b_p:=a_{p+in}a_{p+(i+1)n}$ for all integer $i$).
Then the subgroup generated by the $a_{p,q}$'s and the $b_p$'s is
isomorphic to $B(e,1,n)$. The sublattice $NCP(e,1,n)\hookrightarrow
NCP(e,e,n+1)$ is the lattice of simple elements of a Garside structure
for $B(e,1,n)$, whose atoms are  the  $a_{p,q}$'s and the $b_p$'s.
These atoms generate a monoid $M(e,1,n)$. Actually, up to isomorphism,
the structures of $NCP(e,1,n)$, $B(e,1,n)$ and $M(e,1,n)$ do not
depend on $e\geq 2$. When $e=2$, we recover the dual braid monoid
of type $B_n$.

\section{Zeta polynomials, reflection degrees and Catalan numbers}
\label{section8}
This section contains some complements about the combinatorics of
the lattice $NCP(e,e,n+1)$, which is the lattice of simple elements
in our Garside structure for $B(e,e,n+1)$.

By a \emph{chain of length $N$} in a poset $(P,\leq)$, we mean a finite
weakly increasing
sequence $p_1\leq p_2 \leq p_3 \leq \dots \leq p_N$ of elements of $P$.
When the poset is finite, one may consider the \emph{Zeta function} $Z_P$,
whose value at a positive
integer $N$ is the number of chains of length $N-1$ in $P$.
For example,
$$Z_P(1)=1, \; Z_P(2)=|P|.$$

Let $W$ be an irreducible real reflection group of rank $r$, with 
reflection degrees $d_1\leq d_2 \leq \dots \leq d_r$ (these are the
degrees of homogeneous algebraically independent generators of the algebra
of $W$-invariant polynomial functions; $d_r$, often denoted by $h$, is the
Coxeter number).
Let $P(W)$ be the lattice of simple elements in the dual braid monoid
attached to $W$.
The following general formula was suggested by Chapoton, and has now been proved
using the work of Reiner and Athanasiadis-Reiner for the classical 
types, and computer checks by Reiner for the exceptional types (\cite{chapoton},
\cite{reiner}, \cite{ar}):
$$Z_{P(W)}(X) = \prod_{i=1}^n \frac{d_i + d_n (X-1)}{d_i}.$$
In particular, one has
$$|P(W)|= \prod_{i=1}^n \frac{d_i + d_n }{d_i}.$$
Inspired by the type $A$ situation, the number
$\prod_{i=1}^n \frac{d_i + d_n }{d_i}$ is called \emph{Catalan number}
attached to $W$.

It is natural to expect Chapoton's formula to continue to hold
for our ``dual monoid'' of type $G(e,e,n+1)$.

We set $$Z(1,1,n)(X):=\prod_{k=2}^n\frac{k+n(X-1)}{k},$$
$$Z(e,1,n)(X):=\prod_{k=1}^n\frac{ek+en(X-1)}{ek}= \prod_{k=1}^n\frac{k+
n(X-1)}{k}$$
and 
$$Z(e,e,n+1)(X):=\frac{n+1+en(X-1)}{n+1}\prod_{k=1}^n\frac{k+n(X-1)}{k}.$$
Since $2,3,\dots,n-1,n$ (resp. $e,2e,\dots,e(n-1),en$, resp.
$e,2e,\dots,e(n-1),en,n+1$) are the reflection degrees of $G(1,1,n)$ in its
irreducible reflection representation (resp. $G(e,1,n)$, resp. $G(e,e,n+1)$),
these terms are the right-hand sides in the corresponding Chapoton's
formulae.

One observes the relations
$$Z(e,e,n+1) = (1+\frac{en}{n+1} (X-1)) Z(e,1,n)$$
and
$$Z(e,1,n) = (1+n(X-1))Z(1,1,n).$$

\begin{theorem}
For any $e\geq 2$, $n\geq 1$, we have
$Z_{NCP(1,1,n)}=Z(1,1,n)$, $Z_{NCP(e,1,n)}=Z(e,1,n)$ and 
$Z_{NCP(e,e,n+1)}=Z(e,e,n+1)$.
\end{theorem}

\begin{corollary}
The cardinal of $NCP(e,1,n)$, resp. $NCP(e,e,n+1)$,
is the ``Catalan number'' attached to $G(e,1,n)$, resp. 
$G(e,e,n+1)$.
\end{corollary}

Another useful consequence of the theorem is that the cardinal
of $\Red_T(c)$ may be computed, since it coincides with the number
of strict $n+1$ chains (knowing the number of weak $k$-chains for 
$k=1,\dots,N$, one obtains the number of strict $N$-chains using
a straightforward inversion formula).

Since $G(1,1,n)$, $G(2,1,n)$, $G(2,2,n+1)$ and $G(e,e,2)$ are real,
the corresponding cases are already known (see \cite{chapoton};
the $G(2,2,n+1)=W(D_{n+1})$
case was first proved by Athanasiadis-Reiner, \cite{ar}).
Note that $Z(e,1,n)$ does not depend on $e\geq 2$ -- the poset
$NCP(e,1,n)$ is actually isomorphic to $NCP(2,1,n)$, but is worth 
considering due to its natural relationships with $NCP(e,e,n+1)$, 
which does depend on $e$.

The rest of the section is devoted to the proof of the theorem.
The general strategy is to prove bijectively, at the level of chains, the
relations $Z_{NCP(e,1,n)} = (1+n(X-1))Z_{NCP(1,1,n)}$ and 
$Z_{NCP(e,e,n+1)} = (1+\frac{n+1}{en} (X-1)) Z_{NCP(e,1,n)}$.

Consider the poset morphism $E:NCP(e,1,n)\rightarrow NCP(1,1,n)$ defined
as follows:
for any $u$ in $NCP(e,1,n)$, $E(u)$ is the element of $NCP(1,1,n)$
whose parts are the images of the parts of $u$ by the map $z\mapsto z^e$
(one checks that this
indeed makes sense).
This map $E$ naturally generalises the map
$\Abs:NCP(2,1,n)\rightarrow NCP(1,1,n)$ constructed
by Biane-Goodman-Nica (\cite{bgn}).

For all $\zeta\in\mu_{en}$, we denote by $s_{\zeta}$ the short element
of $NCP(e,1,n)$ whose parts are $\{0\}$, 
$\{\zeta,\zeta_{en}\zeta,\dots,\zeta_{en}^{n-1}\zeta\}$
and their images under the $\mu_e$-action.
The element $s_{\zeta}$ is maximal among short symmetric elements.
There are $n$ such elements (we have $s_{\zeta}=s_{\zeta_e\zeta}$).
The element $l_{\zeta}:=\overline{s_{\zeta}}$ is a minimal long asymmetric
element.

\begin{lemma}
\label{lemmaE}
\begin{itemize}
\item[(i)] The map $E$ is a poset morphism. 
For all $u\preccurlyeq v$ in $NCP(e,1,n)$, one has
$E(u \backslash v) = E(u) \backslash E(v)$. In particular,
$\overline{E(u)}=E(\overline{u})$.
\item[(ii)]
Let $\zeta\in\mu_{en}$.
The map $E$ restricts to a poset isomorphism from $[\disc,s_{\zeta}]$ to
$NCP(1,1,n)$.
\item[(iii)]
Let $\zeta\in\mu_{en}$.
The map $E$ restricts to a poset isomorphism from $[l_{\zeta},\triv]$ to
$NCP(1,1,n)$.
\end{itemize}
\end{lemma}

\begin{proof}
(i) is easy.

(ii).
Consider the restriction map $\res$ from
$[\disc,s_{\zeta}]$ to $NCP_{\{\zeta,\zeta_{en}\zeta,\dots,\zeta_{en}^{n-1}
\zeta\}}$. Clearly, $\res$ is an isomorphism, and the natural map
$NCP_{\{\zeta,\zeta_{en}\zeta,\dots,\zeta_{en}^{n-1}
\zeta\}} \rightarrow NCP(1,1,n)$ is an isomorphism. By composing them,
one obtains the desired isomorphism from
$[\disc,s_{\zeta}]$ to $NCP(1,1,n)$.

(iii) follows from (i) $+$ (ii).
\end{proof}

One could define a category of \emph{complemented lattices},
axiomatising the properties of our lattices and of the operations
$\backslash$. The statement (i) in the above lemma would then express
that $E$ is a \emph{morphism of complemented lattices}. Similarly,
the case (A) of the definition of the complement operation in $NCP(e,e,n+1)$
is designed to make $*$ a morphism of complemented lattices.
The definition below could also be generalised
to study chains in any complemented lattices. 
This axiomatic approach is not necessary to the present work -- the systematic
investigation could however be interesting.

\begin{definition}
Let $P$ be one of the posets $NCP(1,1,n)$, $NCP(e,1,n)$ or $NCP(e,e,n+1)$,
endowed with its operations $\vee$, $\backslash$, etc...
Let $N\in \BZ_{\geq 0}$.
An \emph{$(N+1)$-derived sequence in $P$} is a sequence
$p=(p_0,\dots,p_N)$ of elements of $P$ such that, for all $i$,
$$p_i = (p_1\vee \dots \vee p_{i-1}) \backslash (p_1\vee \dots \vee p_{i})$$
and $p_1\vee \dots \vee p_N=\triv$.

To any $N$-chain $u=(u_1,\dots,u_N)$, we associate an $(N+1)$-derived
sequence $\partial u$, defined by 
$$\partial u := (\disc \backslash u_1=u_1,
u_1\backslash u_2,u_2 \backslash u_3, \dots,u_{N-1}\backslash u_N,
u_N\backslash \triv = \overline{u_N}).$$

To any $(N+1)$-derived sequence $p=(p_0,\dots,p_N)$, we associate
an $N$-chain $\int p$, defined by
$$\int p:= (p_0,p_0\vee p_1,p_0\vee p_1\vee p_2,
\dots,p_0\vee p_1\vee \dots \vee p_{N-1}).$$
\end{definition}

Checking that $\partial u$ is indeed a derived sequence, and that
$\int p$ is indeed a chain, is trivial. It is also trivial to check that

\begin{lemma}
The maps $\partial$ and $\int$ are reciprocal bijections between
the sets of $N$-chains and the set of $(N+1)$-derived sequences in the
corresponding lattice.
We have $\partial E = E \partial$, $\int E = E \int$,
$\partial * = * \partial$, $\int * = * \int$.
\end{lemma}
where we still denote by $E$ the operation sending a sequence 
$(s_1,\dots,s_k)$ to $(E(s_1),\dots,E(s_k))$, and similarly for $*$.

The first section of
\cite{bgn} inspired both the following lemma and its proof.

\begin{lemma}
For any $e,n,N$, the map $E$ from $(N+1)$-derived sequences
in $NCP(e,1,n)$ to $(N+1)$-derived sequences in $NCP(1,1,n)$
is $(1+nN)$-to-$1$ fibration 
(\ie, the pre-image of any $(N+1)$-derived sequence in $NCP(1,1,n)$
has cardinal $1+nN$).
\end{lemma}

\begin{proof}
Observe first that any derived sequence $p=(p_0,\dots,p_N)$
in $NCP(e,1,n)$ contains exactly
one long element: let $i$ be the first integer such
$p_0\vee \dots \vee p_i$ is long; then $p_i$, being the complement
of a short element in a long element, is long. All other elements are
either complements of short elements in short elements, or complements
of long elements in long elements -- in both cases, they must be short
(Lemma \ref{viceversa}).

Let $q:=E(p)$. For all $j$, let $m_j$ be the number of parts of $q_j$.
We have
\begin{equation*}
\begin{split}
\sum_{j=0}^N m_j = \sum_{j=0}^N (n-ht(q_j)) & = (n+1)N -
\sum_{j=0}^N ht(q_j)  \\
& = (n+1)N - ht(q_0\vee \dots \vee q_N) = (n+1)N -ht(\triv)=
nN+1.
\end{split}
\end{equation*}
(the relation $\sum_{j=0}^k ht(q_j) = ht(q_0 \vee \dots \vee q_k)$ is proved
by induction on $k$, using $q_k= (q_0 \vee \dots \vee q_{k-1}) \backslash 
(q_0 \vee \dots \vee q_k)$ and $ht(u) + ht(u\backslash v) = ht(v)$).
Consider the disjoint union of all the parts of all the terms of $q$.
Among these $nN+1$ parts is the image of the long part of $p_i$.

The lemma will be proved if we establish that, choosing any of
$nN+1$ parts (\ie, choosing a given $i$ in $\{0,\dots,N\}$ and
choosing a part of $q_i$)
one may uniquely reconstruct $p\in E^{-1}(q)$ such that $p_i$ is long, 
with long part sent by $E$ to the chosen part.

Suppose we have chosen $q_i$ and a part $a$ of $q_i$.
Let $\zeta\in a$. We want $p_i$ to have a long part containing
the $e$ points in $E^{-1}(\zeta)$. This characterises
a unique element of $p_i\in E^{-1}(q_i)$ (Lemma \ref{lemmaE} (iii)).
We now have to see that the remaining $p_j$'s are then
uniquely determined.
Since $u_i:=p_1\vee \dots \vee p_i$ must 
satisfy $E(u_i)=q_1\vee \dots \vee q_i$ and $p_i\preccurlyeq u_i$, it
is uniquely determined by Lemma \ref{lemmaE} (iii).
Among the constraints is that 
$p_i=(p_1\vee \dots \vee p_{i-1}) \backslash
u_i$. 
It is not difficult to see that, in $NCP(e,1,n)$, for all $u$,
the map $[\disc,u]\rightarrow [\disc,u],x\mapsto x\backslash v$ is a bijection
(this follows from the similar statement in $NCP(1,1,en)$,
which itself easily reduces to the case when $u=\triv$, which itself
follows from the fact that the operation
$x\mapsto \overline{\overline{x}}$ is bijective, which is contained in
Lemma \ref{barbar}).
Applying this to $u_i$,
we find the existence of $x\in NCP(e,1,n)$, uniquely determined,
such that
$p_1\vee \dots \vee p_{i-1} = x\backslash u_i$.
Since $p_i$ is long, so is $u_i$, thus
$p_1\vee \dots \vee p_{i-1}$ is short.
All $p_1,\dots,p_{i-1}$ must be finer than this short element. By Lemma
\ref{lemmaE} (ii), the chain $\int (q_0,\dots,q_{i-1})$ lifts via $E$ to a
unique chain in the interval $[0,x\backslash u_i]$.
This implies that one may uniquely reconstruct $p_0,\dots,p_{i-1}$.

To reconstruct $p_{i+1},\dots,p_N$, similarly observe that
the problem amounts to finding the pre-image of the chain
$(q_1\vee \dots \vee q_i,q_1\vee \dots \vee q_{i+1},\dots,
q_1\vee \dots \vee q_{N+1})$, with the condition that the initial
term is lifted to $u_i$. Since $u_i$ is long, this problem
admits a unique solution (Lemma \ref{lemmaE} (iii)).
\end{proof}

The $NCP(e,1,n)$ case of the theorem follows from the above lemma
(rephrased to deal with chain rather than derived sequences) and the
relation $Z(e,1,n) = (1+n(X-1))Z(1,1,n)$.

We now deal with the $NCP(e,e,n+1)$ case.
There are two sorts of chains in $NCP(e,e,n+1)$:
\emph{symmetric chains}
(chains consisting only of symmetric partitions) and
\emph{asymmetric chains} (chains containing at least an asymmetric 
partition).
Via $*$, symmetric $N$-chains in $NCP(e,e,n+1)$
are in $1$-to-$1$ correspondence with
$N$-chains in $NCP(e,1,n)$.
Expanding the relation
$Z(e,e,n+1) = (1+\frac{en}{n+1} (X-1)) Z(e,1,n)$
and using the $NCP(e,1,n)$ case already proved, we are
left with having to prove that the number of 
asymmetric $N$-chains is $\frac{en}{n+1} N  Z(e,1,n) =
\frac{en}{n+1} N Z(2,1,n)$.
Using the $NCP(2,2,n+1)$ case already proved by Athanasiadis-Reiner,
\cite{ar},
the claimed result follows from:

\begin{lemma}
The number of asymmetric $N$-chains in $NCP(e,e,n+1)$ is
$e/2$ times the number of asymmetric $N$-chains in $NCP(2,2,n+1)$.
\end{lemma}

\begin{proof}
Let $u$ be an asymmetric partition, let $\zeta$ be the successor
of $0$ for the counter-clockwise cyclic ordering of the asymmetric part.
We say that $\zeta$ is the \emph{direction of $u$}.
The \emph{direction} of an asymmetric chain is the direction of the first
asymmetric partitions appearing in the chain.
Clearly, there are $ne$ possible directions for asymmetric $N$-chains,
and these chains are equidistributed according to the possible directions.
To prove the lemma, it is enough to prove that the number
of asymmetric $N$-chains in $NCP(e,e,n+1)$ with direction $1$ depends
only on $n$ and $N$, and not on $e$.

To prove this, observe that
an asymmetric chain with direction $1$ is uniquely determined
by its restriction to $\{0,\zeta_{en}^{-e+1},\zeta_{en}^{-e+2},
\dots,1,\dots,\zeta_{en}^e\}$ and 
use the map ``$E^{-1}\circ E$'' from $NCP(e,e,n+1)$ to $NCP(2,2,n+1)$.
\end{proof}

\section{A remark on Brou\'e-Malle-Rouquier generators}
\label{section9}

In this section, we will show that the BMR presentation defines a monoid which
does not embed in the group defined by the same presentation,
and that indeed, the
submonoid of the braid group of $G(e,e,n+1)$ generated by the BMR generators is
not finitely presented.

Let $\cal T$ denote the BMR generators, and $\cal R$ the BMR relations.
For ease on the eye, we will write $1$ for $\tau_2'$ ($\leftrightarrow a_n$), 
$2$ for $\tau_2$ ($\leftrightarrow a_0$)  and
$i$ for $\tau_i$ ($\leftrightarrow a_{i-3,i-2}$) for $i=3,\ldots, n-1$.
We now deduce some relations in $B(e,e,n+1)$. Firstly
$$
\begin{array}{rcl}
2\ 1  3213 \langle 21 \rangle^{e-2} 
        & = & 321321 \langle 21 \rangle^{e-2} \ = \ 3213 \langle 21 \rangle^{e} 
         \ = \ 3213 \langle 12 \rangle^{e} \\
         & = & 3213 \ 1 \langle 21 \rangle^{e-1} \ = \ 32 313 \langle 21 \rangle^{e-1} 
         \ =  \ 2 \ 3213 \langle 21 \rangle^{e-1} ,
\end{array}
$$
so by left cancellation, $1 3 2 1 3 \langle 2 1 \rangle^{e-2}  = 3 2 1 3 \langle 2 1 \rangle^{e-1}$.
This new relation can be written $1 w = w x $ where $w$ is $3213 \langle 21 \rangle^{e-2}$ 
and $x$ is $2$ if $e$ is even, and $1$ if $e$ is odd. For the same letter $x$,
$\langle 21 \rangle^e x \equiv \langle 21 \rangle^{e+1} \equiv 2 \langle 12 \rangle^e$
hence for any $k$,
$$21^k 3213\langle 21 \rangle^{e-2} \equiv 21^k w 
= 21 w x^{k-1}  = 3213 \langle 21 \rangle^e x^{k-1}
= 3213 2^{k-1} \langle 21 \rangle^e  \equiv  3213 2^k 1 \langle 21 \rangle^{e-2},$$
so by right cancellation this time, $2 1^k 3 2 1 3  = 3 2 1 3 2^k 1$.


Let $M_{\cal T}$ denote the submonoid of the braid group of $G(e,e,n+1)$ generated by the BMR generators.
Then $M_{\cal T}$ is a submonoid of $M=M(e,e,n+1)$, the monoid defined by the presentation given at the
beginning of this section. We have a solution to the word problem for $M$ and can calculate least common multiples, 
and from this we will be able to deduce that there can be no finite presentation for $M_{\cal T}$.
Suppose there were such a finite presentation; then there would be an $l$
for which all relations -- which must be words over $\cal T$ -- are of length at most $l$. 
For all $k$, we have
$21^k3213 = 32132^k1$, so a relation must be applicable to the
word $32132^k1$. Since $k$ can be arbitrarily large, this means that a relation
must be applicable to either $32132^k$ or to $2^k 1$. That is, there would have to exist some
(different) word $W$ over $\cal T$ such that $32132^k = W$ or $2^k 1 = W$. We will see that this
cannot be the case.

%
%

%

There is a surjection $B(e,e,n+1) \rightarrow B_n$, the braid group of
type $A$ on $n$ strings, given by $1,2 \mapsto \sigma_1$ and
 $i \mapsto \sigma_{i-1}$ for $i \geq 3$. Thus 
$2^k 1\stackrel{\varphi}{\mapsto} \sigma_1^{k+1}$; so $2^k 1$ 
can only be rewritten in terms of $\{1,2\}$. Suppose that 
$2^k 1 = w 2$ for some word $w$ over $\{2,1\}$; then
$2^k 1 \equiv 2^{k-1} 21 = 2^{k-1}  a_{-n} 2 $; so by right cancellation,
$w = 2^{k-1} a_{-n}$. However 
$2^{k-1} a_{-n}  \equiv a_{0}^{k-1} a_{-n}$
is in a singleton equivalence class in $M$, so can never be 
rewritten in terms of $\{1,2\}$. Thus it remains to 
show that $32132^k$ cannot be rewritten in terms
of $\cal T$. 

Let $U$ denote the word $3213\,2^k$. It suffices to show that 
we cannot rewrite $U$ in terms of ${\cal T}_3 := \{1,2,3\}$.
Since $M$ is Garside, it has the `reduction property', and
so we can quickly determine divisibility in $M$ by a generator using 
the method of $a$-chains (see~\cite{corran}). 

Observe that $U = u \, 3213$ where $u = a_{-n}^{k}$; 
$u$ is not divisible by $1$, $2$, or $3$, so is certainly not 
rewriteable over these letters. Also, 
neither $1$ nor $2$ right divides
$u \, 3$ so $U$ cannot be ${\cal T}_3$-rewritten to end with $213$.
Since $3 2 = a_1 3$, we have $U = v \,313$ where $v = u a_{1}$,
which is in a singleton equivalence class in $M$, so not ${\cal T}_3$-rewriteable; 
so $U$ cannot be ${\cal T}_3$-rewritten to end with $313$.
Furthermore, $v\,3$ is not right divisible by $1$, so  $U$ cannot be ${\cal T}_3$-rewritten to end with $113$.
Thus $U$ cannot be ${\cal T}_3$-rewritten to end with $13$.

Next, write $U =  w \, 323 $ where $w = u a_{1-n}$; 
$w$ is not  ${\cal T}_3$-rewriteable; and $w \, 3$ is not right divisible by 1 or 2. 
Thus $U$ cannot be rewritten over $\{1,2,3\}$ to end with  $23$.
Finally, $U = x  \,3^2 $ where $x = u  a_{1-n} a_1$;
$x$ is not right divisible by 1,2 or 3, so we have that $U$ cannot be rewritten
to end with $33$ either. This exhausts all the possibilites for ${\cal T}_3$-rewriting
$U$ to end with $3$.

Now write $U = y  \,321$ where $y= u a_{1,n}$. This is not
right divisible by $1,2$ or $3$; so not ${\cal T}_3$-rewriteable over these letters. 
Moreover, $y \,3$ is not right divisible by $1$ or $2$; and $y \,32$ is not right divisible by $1$.
Thus we have that $U$ cannot be rewritten to end with $12$ or $11$. 
Write $U = v \, 131$ where $v = u a_{1}$, 
which is not right divisible by $1, 2$ or $3$; furthermore, $v1$ is not divisible
by $2$ or $3$. In this way we have exhausted all the possibilites for rewriting
$U$ over $\{1,2,3\}$ ending with $1$.

Finally, we want to show that $U$ cannot be rewritten to end with $2$. By 
right cancellation, assume that $k=0$. Now $3213 = a_{1-n} 232$, but 
$a_{1-n} 23$ is not left divisible by $1, 2$ or $3$,
so cannot be rewritten over $\{1,2,3\}$. Thus

\begin{proposition}
There is no finite presentation for the submonoid of the braid group of type $G(e,e,r)$
generated by the Brou\'e-Malle-Rouquier generators.
\qed
\end{proposition}

\end{document}